\documentclass[10pt, a4paper]{amsart}
\usepackage{amssymb,graphicx}
\usepackage{mathrsfs, bbm, slashed}
\usepackage{color}
\usepackage{enumerate}
\usepackage{a4wide}

\usepackage[utf8]{inputenc}
\usepackage[colorlinks=true,citecolor=blue,linkcolor=blue]{hyperref}

\newcommand{\R}{\mathbb{R}}
\newcommand{\C}{\mathbb{C}}
\newcommand{\Z}{\mathbb{Z}}
\newcommand{\N}{\mathbb{N}}

\newcommand{\sA}{\mathscr{A}}
\newcommand{\sB}{\mathscr{B}}
\newcommand{\sE}{\mathscr{E}}
\newcommand{\sL}{\mathscr{L}}
\newcommand{\sH}{\mathscr{H}}
\newcommand{\sK}{\mathscr{K}}
\newcommand{\sS}{\mathscr{S}}
\newcommand{\sW}{\mathscr{W}}

\newcommand{\fp}{\mathfrak{p}}

 \newcommand{\psidos}{$\Psi$DOs}

\newcommand{\op}{\operatorname} 

\newcommand{\Sp}{\op{Sp}}
\newcommand{\Tr}{\op{Tr}}
\newcommand{\dist}{\op{dist}}
\newcommand{\ran}{\op{ran}}

\newcommand{\Hol}{\op{Hol}}

\newcommand{\car}{\mathbbm{1}}

\newcommand{\scal}[2]{\ensuremath{\left\langle #1 | #2 \right\rangle}} 
\newcommand{\bigscal}[2]{\ensuremath{\big\langle #1 | #2 \big\rangle}} 

\newcommand{\dom}{\op{dom}}

\numberwithin{equation}{section}

\newtheorem{theorem}{Theorem}[section]
\newtheorem{proposition}[theorem]{Proposition}
\newtheorem*{conditionW}{Condition (W)}

\newtheorem{lemma}[theorem]{Lemma}

\newtheorem{theoremalph}{Theorem}
\newtheorem{propositionalph}[theoremalph]{Proposition}
\newtheorem{corollaryalph}[theoremalph]{Corollary}

\theoremstyle{definition}
\newtheorem{definition}[theorem]{Definition}
\newtheorem{condition}[theorem]{Condition}

\theoremstyle{remark}
\newtheorem{remark}[theorem]{Remark}

\newcommand{\shD}{\slashed{D}}
\newcommand{\shS}{\slashed{S}}

\title{Noncommutative Geometry, Spectral Asymptotics, and Semiclassical Analysis, II. Operator-Theoretic Tools}
\date{}

\author{Rapha\"el Ponge}
 \address{Department of Mathematics and Statistics, University of Wyoming, Laramie, WY, USA}
 \address{Department of Mathematics and Statistics, University of Ottawa, Canada}
 \email{ponge.math@icloud.com}

\begin{document}
\begin{abstract}
This is the second part of a two-part work on semiclassical Weyl laws and Connes' integration formulas for spectral triples. It establishes three operator-theoretic results, which are used in Part~I but are of independent interest. The first is a factorization formula for the fractional commutators $[|D|^{-\alpha},a]$, which is obtained by means of double operator integrals. The second is an enhancement, for arbitrary unital spectral triples, of a lemma of Connes--Sukochev--Zanin on the differences $(a^{1/2}|D|^{-\alpha}a^{1/2})^s-|D|^{-\alpha s}a^s$, together with a holomorphic version of it. The third is a Tauberian theorem showing that zeta function and heat trace conditions of the type sought by Connes imply Weyl laws for the operators $aD^2a$, where $a$ ranges over positive invertible elements of the algebra. Along the way, we obtain a version for spectral triples of the Tauberian theorem of McDonald--Sukochev--Zanin without any extra regularity assumption or restriction on the degree of summability. This paper can be read largely independently of Part~I.
\end{abstract}

\maketitle

\tableofcontents

\section{Introduction}\label{sec:Intro}
This paper is the second part of a two-part work on semiclassical Weyl laws and Connes' integration formulas for spectral triples. In Part~I~\cite{Po:Part1} these results are derived from a spectral condition, Condition~(W), which requires Weyl laws for the operators $aD^2a$, where $a$ ranges over positive invertible elements of the algebra of the spectral triple. The approach of Part~I relies on three operator-theoretic results: a factorization lemma for fractional commutators, an enhancement of a lemma of Connes--Sukochev--Zanin~\cite{CSZ:MS17}, and a Tauberian theorem. The aim of this paper is to establish these results. They are respectively the Factorization Lemma (Theorem~\ref{thm:Intro.factorization}), the Enhanced CSZ Lemma (Theorem~\ref{thm:Intro.CSZ}) and the Tauberian Theorem (Theorem~\ref{thm:Intro.Tauberian}) below. We also establish a holomorphic version of the Enhanced CSZ Lemma (Theorem~\ref{thm:Intro.CSZhol}), which is a key ingredient in the proof of the Tauberian Theorem. These results are of independent interest, and this paper can be read largely independently of Part~I. The only results of Part~I that are used are a few auxiliary results, which are recalled in Section~\ref{chap:Tauberian}.

Throughout this paper we let $(\sA,\sH,D)$ be a $p$-summable spectral triple, where $p$ may be \emph{any} positive number (see Section~\ref{sec:Prelim} for the precise definitions). We denote by $\Pi_0$ the orthogonal projection onto $\ker D$, by $|D|^{-1}$ the partial inverse of $|D|$, and by $F=D|D|^{-1}$ the sign of $D$. We also denote by $\overline{\sA}$ the closure of $\sA$ in $\sL(\sH)$, and we let $\sA_{++}$ be the set of elements of $\sA$ that are positive and invertible in $\overline{\sA}$. In addition, if $X$ is a positive compact operator, we denote by $\lambda_0(X)\geq \lambda_1(X)\geq \cdots$ its eigenvalues counted with multiplicity.

\subsection{Factorization Lemma}
The first main result of this paper is the following factorization formula for the fractional commutators $[|D|^{-\alpha},a]$, $a\in \sA$.

\begin{theoremalph}[Factorization Lemma]\label{thm:Intro.factorization}
 Let $\alpha, \beta, \gamma\in [0,1]$ be such that, either $\beta+\gamma<\alpha+1$, or $\beta+\gamma=\alpha+1$ and $\max(\beta,\gamma)<1$.
 \begin{enumerate}
 \item[(i)] There is a continuous linear operator $\Phi_{\beta,\gamma}^{\alpha}:\sL(\sH)\rightarrow \sL(\sH)$ such that, for all $a\in \sA$, we have
 \begin{equation}\label{eq:Intro.factorization}
  \left[ |D|^{-\alpha},a\right]  =|D|^{-\beta}\Phi_{\beta,\gamma}^{\alpha}\left([D,a]\right)|D|^{-\gamma}
  + |D|^{-(\alpha+1)}F[D,a]\Pi_0 +\Pi_0[D,a]F|D|^{-(\alpha+1)}.
\end{equation}

\item[(ii)] Assume further that $\beta+\gamma<\alpha+1$. If $\max(\beta,\gamma)<1$, then $\Phi_{\beta,\gamma}^{\alpha}$ gives rise to a continuous linear operator,
\begin{equation*}
 \Phi_{\beta,\gamma}^{\alpha}: \sL(\sH)\longrightarrow \sL_{r,\infty}, \qquad \text{with}\ r:=(1+\alpha -\beta-\gamma)^{-1}p.
\end{equation*}
If $\max(\beta,\gamma)=1$, then the same result holds for any $r>(1+\alpha -\beta-\gamma)^{-1}p$.
\end{enumerate}
\end{theoremalph}

The above result seems to be new. However, its proof relies on basic techniques of double operator integrals (DOIs) introduced by Birman--Solomyak~\cite{BS:PMF66, BS:PMF73, BS:IEOT03} in the 60s (see also~\cite{ST:LNM19}). The possibility of factoring out $|D|^{-\gamma}$ with $\gamma=\alpha$, including for $\alpha=1$, is crucial in the proof of the Enhanced CSZ Lemma below (see Remark~\ref{rmk:CSZ.alpha}).

As an application of Theorem~\ref{thm:Intro.factorization}, the following commutator estimates are established in Part~I. 

\begin{propositionalph}[Commutator Estimates; see Part~I~{\cite[Proposition~I]{Po:Part1}}]\label{prop:Intro.Com-Dq}
Given any $q>0$, for all $a\in \sA$, we have
 \begin{equation*}
 \left[ |D|^{-q},a\right] \in \sL_{(q+1)^{-1}p,\infty}.
\end{equation*}
In fact, there is a constant $C_{D,q}>0$, which depends only on $D$ and $q$, such that
\begin{equation*}
 \left\| \left[ |D|^{-q},a\right]\right\|_{(q+1)^{-1}p,\infty}\leq C_{D,q}  \left\| \left[ D,a\right]\right\| \qquad \forall a \in \sA.   
\end{equation*} 
\end{propositionalph}

Proposition~\ref{prop:Intro.Com-Dq} subsumes and improves various commutator estimates for spectral triples (see, e.g., \cite{BJ:CRAS83, CM:CRAS86, CP:CJM98, SWW:CRAS98}), and holds, e.g., for all Lipschitz functions in the case of Dirac spectral triples.

\subsection{Enhanced CSZ Lemma}
The second main result is an enhancement for (unital) spectral triples of a deep lemma of Connes--Sukochev--Zanin~\cite[Lemma~5.3]{CSZ:MS17} (hereafter called the CSZ Lemma). The original CSZ Lemma, stated without proof by Connes~\cite[Lemma~IV.3.11]{Co:NCG} and proved in~\cite{CSZ:MS17}, asserts in our setting that, for $0<q<p$, the differences $a^{p/2}|D|^{-p}a^{p/2}-(a^{1/2}|D|^{-q}a^{1/2})^{p/q}$, $a\in \sA_{++}$, are in $(\sL_{1,\infty})_0$. Our approach relies on an integral representation for differences of the form $(A^{1/2}BA^{1/2})^z-B^zA^z$, $\Re z>1$, which was established by Connes--Sukochev--Zanin~\cite[Lemma~5.2]{CSZ:MS17} and further refined by Sukochev--Zanin~\cite[Theorem~5.2.1]{SZ:Ast23}. Combining it with the Factorization Lemma, we obtain the following enhancement of the CSZ Lemma.

\begin{theoremalph}[Enhanced CSZ Lemma]\label{thm:Intro.CSZ}
 Given any $\alpha \in (0,1]$ and $s>1$, for every $a\in {\sA}_{++}$, we have
 \begin{equation*}
 \big( a^{\frac12}|D|^{-\alpha} a^{\frac12}\big)^{s} - |D|^{-\alpha s} a^{s} \in \sL_{(\alpha s +1-\epsilon)^{-1}p,\infty} \qquad \forall \epsilon\in (0,1).
\end{equation*}
\end{theoremalph}

Each of the operators $(a^{1/2}|D|^{-\alpha}a^{1/2})^{s}$ and $|D|^{-\alpha s}a^{s}$ is contained in $\sL_{(\alpha s)^{-1}p,\infty}$. Theorem~\ref{thm:Intro.CSZ} thus shows that their difference lies in the much smaller ideals $\sL_{(\alpha s+1-\epsilon)^{-1}p,\infty}$, $\epsilon\in(0,1)$. In particular, this difference is trace-class as soon as $\alpha s>p-1$. We also obtain the following holomorphic version of the Enhanced CSZ Lemma.

\begin{theoremalph}[Holomorphic Enhanced CSZ Lemma]\label{thm:Intro.CSZhol}
 Let $a\in \sA_{++}$ and $\alpha\in (0,1]$, and set
 \begin{equation*}
 T(s):=\big(a^{\frac12}|D|^{-\alpha} a^{\frac12}\big)^{s} - |D|^{-\alpha s} a^{s}, \qquad \Re s>1.
\end{equation*}
 Given any $q>0$, set $c:=\max(1,\alpha^{-1}(q^{-1}p-1))$. Then $T(s)$, $\Re s>c$, is a holomorphic family in $\sL_q$, i.e., $T(s)\in \Hol\left(\Re s>c;\sL_q\right)$.
\end{theoremalph}

In the setting of unital spectral triples, a version of Theorem~\ref{thm:Intro.CSZhol} for $\alpha=1$ and $q=1$ is given by Sukochev--Zanin~\cite[Corollary~5.4.11]{SZ:Ast23}, where it is shown that $(A^{1/2}BA^{1/2})^s-B^sA^s\in\sL_1$ for $\Re s>p-1$, under strong regularity assumptions and the summability condition $p>2$. Theorem~\ref{thm:Intro.CSZhol} holds without these assumptions, and for all $\alpha\in(0,1]$. 

Theorem~\ref{thm:Intro.CSZhol} has the following consequence for zeta functions, which is a key ingredient in the proof of the Tauberian theorem below.

\begin{corollaryalph}\label{cor:Intro.zeta}
 Let $a\in \sA_{++}$ and $\alpha\in (0,1]$ be such that $\alpha<p$, and set $\delta:=\min(1,p-\alpha)$. Then the function
 \begin{equation*}
 \Tr\Big[\big(a^{\frac12}|D|^{-\alpha} a^{\frac12}\big)^{s}\Big]-\Tr\big[a^{s}|D|^{-\alpha s}\big], \qquad \Re s>\alpha^{-1}p,
\end{equation*}
has a holomorphic extension to the half-plane $\Re s>\alpha^{-1}(p-\delta)$.
\end{corollaryalph}

In the setting of unital spectral triples, Corollary~\ref{cor:Intro.zeta} is a refinement of a result of Sukochev--Zanin~\cite[Theorem~5.4.2]{SZ:Ast23}. The latter result is proved for $\alpha=1$ and $p>2$ under some strong regularity assumptions (it actually holds for non-unital spectral triples). A weaker version, which only provides a continuous extension to the closed half-plane $\Re s\geq p$, is given by McDonald--Sukochev--Zanin~\cite[Theorem~2.2]{MSZ:LMP22} for Lipschitz-regular spectral triples with $p>2$. All these regularity conditions are removed in Corollary~\ref{cor:Intro.zeta}.

More importantly, Corollary~\ref{cor:Intro.zeta} does not involve any restriction on the degree of summability, since we can always choose $\alpha\in(0,1]$ small enough so that $\alpha<p$, and hence $\delta=\min(1,p-\alpha)>0$. This is what allows us to remove the condition $p>2$ in the Tauberian theorem below.

\subsection{Tauberian Theorem}
Let $\tau:\overline{\sA}\rightarrow \C$ be a non-zero positive linear map. In Part~I the semiclassical Weyl laws and integration formulas are derived from the following spectral condition.

\begin{conditionW}
 For every $a\in {\sA}_{++}$, we have
\begin{equation}\label{eq:Intro.conf-Weyl}
 \lim_{j\rightarrow \infty}j^{-\frac{2}{p}}\lambda_j\big(aD^2a\big)= \tau\big[a^{-p}\big]^{-\frac{2}{p}},
\end{equation}
where $\lambda_0(aD^2a)\leq \lambda_1(aD^2a)\leq \cdots$ are the positive eigenvalues of $aD^2a$ counted with multiplicity (see Section~\ref{subsec:Prelim.ST}).
\end{conditionW}

McDonald--Sukochev--Zanin~\cite{MSZ:LMP22} derived semiclassical Weyl laws for spectral triples from a Tauberian theorem based on some condition on the zeta functions $\Tr[a^s|D|^{-s}]$, where $a$ ranges over all positive elements of the $C^*$-algebra $\overline{\sA}$ (see Condition~\ref{cnd:TauberianMSZ} below). Their Tauberian theorem further requires Lipschitz regularity and $p$-summability with $p>2$.

During an online seminar in May 2021, Alain Connes pointed out the need to relate this condition to more standard Tauberian conditions, in terms of the zeta functions $\Tr[a|D|^{-s}]$ or the heat traces $\Tr[ae^{-tD^2}]$, $t>0$, which do not involve powers of $a$. He also pointed out that closedness under holomorphic functional calculus should play a role.

Our third main result provides Tauberian conditions of this type, which only need to be checked on a holomorphically closed Fr\'echet subalgebra containing $\sA$, such as $\sA$ itself when it is holomorphically closed. They imply Condition~(W) without any extra regularity or summability assumption.

\begin{theoremalph}[Tauberian Theorem]\label{thm:Intro.Tauberian}
 Assume there exists a holomorphically closed Fr\'echet subalgebra $\sB\subseteq \overline{\sA}$ containing $\sA$ (see Section~\ref{subsec:Prelim.ST}) such that one of the following conditions is satisfied:
\begin{itemize}
\item \textup{\textbf{Condition (Z)}}. For every $a\in \sB$, the function
\begin{equation*}
\Tr\big[a|D|^{-s}\big] -p\tau(a)(s-p)^{-1}, \qquad \Re s>p,
\end{equation*}
 has a continuous extension to the closed half-plane $\Re s\geq p$.\smallskip

  \item \textup{\textbf{Condition (H)}}. There is $\delta>0$ such that, for all $a\in \sB$, as $t\rightarrow 0^+$ we have
 \begin{equation*}
  \Tr\big[ae^{-tD^2}\big]=  \Gamma\left(1+\frac{p}{2}\right) t^{-\frac{p}{2}}\left[\tau(a) + \op{O}\big(t^{\delta}\big)\right].
\end{equation*}
\end{itemize}
Then Condition~\textup{(W)} is satisfied.
\end{theoremalph}

Conditions~(Z) and~(H) are the types of Tauberian conditions sought by Alain Connes. In practice, they are easier to check than the Tauberian condition of~\cite{MSZ:LMP22} and hold in much greater generality. The proof of Theorem~\ref{thm:Intro.Tauberian} shows that each of these conditions implies a weaker version of the Tauberian condition of~\cite{MSZ:LMP22}, Condition~($\textup{Z}_0$), which in turn implies  Condition~(W) (see Section~\ref{chap:Tauberian}). We thus obtain a stronger version for spectral triples of the Tauberian theorem of~\cite{MSZ:LMP22} without its Lipschitz regularity assumption and restriction on the degree of summability (see Remark~\ref{rmk:Tauberian.comparison-MSZ}).
\subsection{Organization of the paper}
The rest of the paper is organized as follows. In Section~\ref{sec:Prelim}, we gather background on weak Schatten classes, complex powers of positive operators, and spectral triples. In Section~\ref{chap:DOIs}, we review the double operator integrals that we need. In Section~\ref{sec:Factorization}, we prove the Factorization Lemma (Theorem~\ref{thm:Intro.factorization}). In Section~\ref{chap:CSZ}, we prove the Enhanced CSZ Lemma (Theorem~\ref{thm:Intro.CSZ}), its holomorphic version (Theorem~\ref{thm:Intro.CSZhol}), and Corollary~\ref{cor:Intro.zeta}. In Section~\ref{chap:Tauberian}, we prove the Tauberian Theorem (Theorem~\ref{thm:Intro.Tauberian}). In Appendix~\ref{app:hol-extension}, we show that the algebra of any spectral triple can be enlarged to a holomorphically closed algebra without changing the operator $D$. In Appendix~\ref{app:BS}, we include a proof of Proposition~\ref{lem:Com.estimate-Delta-weak}, which is due to Birman--Solomyak.

\subsection*{Acknowledgements}
The author wishes to thank Alan Carey, Magnus Goffeng, Edward McDonald, Adam Rennie, Grigori Rozenblum, Fedor Sukochev, and Teun van Nuland for various discussions related to the subject matter of this article. The research for this article was partially funded by NSFC grant No.~11971328 (China).

\section{Preliminaries}\label{sec:Prelim}
In this section, we gather background material on weak Schatten classes, complex powers of positive operators, and spectral triples.

\subsection{Weak Schatten classes}\label{subsec:Prelim.Schatten}
We refer to~\cite{GK:AMS69, Si:AMS05} for more detailed accounts on Schatten classes and weak Schatten classes. Let $\sK$ be the ideal of compact operators on $\sH$. Given any operator $T\in \sK$, we denote by $(\mu_j(T))_{j\geq 0}$ its sequence of \emph{singular values}, i.e., $\mu_j(T)$ is the $(j+1)$-th eigenvalue counted with multiplicity of the absolute value $|T|=\sqrt{T^*T}$. If $T\geq 0$, then $\mu_j(T)=\lambda_j(T)$ for all $j\geq 0$.

By definition the \emph{Schatten class} $\sL_q$, $q>0$, consists of operators $T\in \sK$ such that $|T|^q$ is trace-class. It is equipped with the quasi-norm,
\begin{equation*}
 \|T\|_{q}:=\Tr\big(|T|^q\big)^{\frac1q}= \bigg( \sum_{j\geq 0}\mu_j(T)^q\bigg)^{\frac1q}, \qquad T\in \sL_q.
\end{equation*}
We obtain a quasi-Banach ideal, which is a Banach ideal when $q\geq 1$. In particular, $\sL_1$ is the ideal of trace-class operators.

The \emph{weak Schatten class}  $\sL_{q,\infty}$, $q>0$, is defined by
\begin{equation*}
 \sL_{q,\infty}:=\left\{T\in \sK; \ \mu_j(T)=\op{O}\big(j^{-\frac1q}\big)\right\}.
\end{equation*}
This is a quasi-Banach ideal with respect to the quasi-norm,
\begin{equation*}
 \|T\|_{q,\infty}:=\sup_{j\geq 0}\;(j+1)^{\frac{1}{q}}\mu_j(T), \qquad T\in \sL_{q,\infty}.
\end{equation*}
We denote by $(\sL_{q,\infty})_{0}$ the closure in $\sL_{q,\infty}$ of the ideal of finite-rank operators. Equivalently, $(\sL_{q,\infty})_{0}$ consists of the operators $T\in \sK$ such that $\mu_j(T)=\op{o}(j^{-1/q})$. We have the following continuous inclusions,
\begin{equation*}
 \sL_q \subseteq  \big(\sL_{q,\infty}\big)_{0}\subseteq \sL_{q,\infty}  \subseteq \sL_{q'}, \qquad 0<q<q'.
\end{equation*}

We will use the following version of H\"older's inequality for weak Schatten classes.

\begin{proposition}[H\"older's Inequality; see~\cite{GK:AMS69, Si:AMS05, SZ:PAMS21}]\label{prop:Hoelder}
 Suppose that $q_1^{-1}+q_2^{-1}=r^{-1}$. If $S\in \sL_{q_1,\infty}$ and $T\in \sL_{q_2,\infty}$, then $ST\in \sL_{r,\infty}$, and we have
\begin{equation*}
 \|ST\|_{r,\infty} \leq  C_{q_1q_2}\|S\|_{q_1,\infty}  \|T\|_{q_2,\infty},
\end{equation*}
where the constant $C_{q_1q_2}>0$ depends only on $q_1$ and $q_2$. If in addition $S\in (\sL_{q_1,\infty})_0$ or $T\in (\sL_{q_2,\infty})_0$, then $ST\in (\sL_{r,\infty})_0$.
\end{proposition}

Given $q>0$, we denote by $\sW^+_{q,\infty}$ the set of positive operators $A\in \sL_{q,\infty}$ such that the limit $\Lambda_q(A):=\lim_{j\rightarrow \infty} j^{\frac1q}\lambda_j(A)$ exists. We will use the following perturbation result of Birman--Solomyak.

\begin{proposition}[Birman--Solomyak~\cite{BS:FAA70, BS:Book}; see also~\cite{Fa:PNAS51}]\label{prop:Bir-Sol}
 Let $q>0$.
 \begin{enumerate}
 \item[(i)] $\sW^+_{q,\infty}$ is a closed subset of $\sL_{q,\infty}$, and $\Lambda_q$ is a continuous functional on $\sW^+_{q,\infty}$.
 \item[(ii)] If $A\in \sW^+_{q,\infty}$ and $B\in \sL_{q,\infty}$ is a positive operator such that $A-B\in (\sL_{q,\infty})_0$, then $B\in \sW^+_{q,\infty}$ and $\Lambda_q(B)=\Lambda_q(A)$.
 \end{enumerate}
\end{proposition}

Let $\sE$ be a quasi-Banach space (e.g., $\sL_q$ or $\sL_{q,\infty}$), and let $\Omega$ be an open subset of $\C$. A map $f:\Omega\rightarrow \sE$ is called holomorphic if it is complex differentiable at every point of $\Omega$ with respect to the quasi-norm topology of $\sE$. In this case, we also say that $f(z)$, $z\in \Omega$, is a holomorphic family in $\sE$. We denote by $\Hol(\Omega;\sE)$ the space of holomorphic maps from $\Omega$ to $\sE$.

\subsection{Complex powers of positive operators}\label{subsec:Prelim.powers}
In what follows, if $X$ is a positive compact operator, then we define $X^z$, $\Re z\geq 0$, by standard continuous or Borel functional calculus. Thus, if $(\xi_\ell)_{\ell \geq 0}$ is an orthonormal eigenfamily such that $X\xi_\ell =\lambda_\ell\xi_\ell$ for all $\ell \geq 0$ (where $\lambda_0\geq \lambda_1\geq \cdots$ are the eigenvalues of $X$ counted with multiplicity), then
\begin{equation}\label{eq:Com.Xz-definition}
 X^z\xi_\ell=\lambda_\ell^z \xi_\ell \quad \ell\geq 0, \qquad X^z=0 \ \text{on}\ \ker X. 
\end{equation}

Similarly, if $A$ is a positive operator on $\sH$ with compact resolvent, we define its powers $A^z$, $z\in \C$, by using the Borel functional calculus for $A$ with the convention that $A^z=0$ on $\ker A$. That is, $A^z=f_z(A)$ with $f_z(t)=\car_{[c,\infty)}(t)t^z$, where $0<c<\dist(0,\Sp(A)\setminus 0)$. Equivalently, if $(\xi_\ell)_{\ell\geq 0}$ is any orthonormal eigenbasis of $\sH$ with $A\xi_\ell=\lambda_\ell\xi_\ell$ for all $\ell\geq 0$, then
\begin{equation*}
 A^z\xi_\ell = \left\{ 
\begin{array}{cl}
 \lambda_\ell^z \xi_\ell & \text{if}\ \lambda_\ell>0,\\
 0 & \text{if}\ \lambda_\ell=0. 
\end{array}\right. 
\end{equation*}
For $\Re z\leq 0$ the operator $A^z$ is bounded. For $\Re z>0$ this is an unbounded normal operator whose domain is
\begin{equation*}
 \dom(A^z)= \left\{ \xi\in \sH; \ \sum |\lambda_\ell(A)|^{2\Re z} \left|\scal{\xi_\ell}{\xi}\right|^2<\infty\right\}. 
\end{equation*}
In any case we have, 
\begin{equation*}
 A^{0}=1-\Pi_0(A), \qquad A^{z_1+z_2}=A^{z_1}A^{z_2}, \quad z_j\in \C. 
\end{equation*}
Here $\Pi_0(A)$ is the orthogonal projection onto $\ker A$, and the second equality holds on $\dom(A^{z_2})\cap \dom(A^{z_1+z_2})$. In particular, $A^{-1}$ is the partial inverse of $A$, and in view of~(\ref{eq:Com.Xz-definition}) we have
\begin{equation*}
 A^{-z}=(A^{-1})^z \qquad \text{if}\ \Re z\geq 0. 
\end{equation*}
Note also that with the above convention $A^{it}$, $t\in \R$, is not a unitary operator unless $\ker A=0$.

\begin{lemma}\label{lem:Tauberian.powers}
 Let $X$ be a (non-zero) positive compact operator on $\sH$. The following holds.
 \begin{enumerate}
\item[(i)] The family $X^{z}$, $\Re z>0$, is a holomorphic family of compact operators. 

\item[(ii)] Assume further that $X\in \sL_{p,\infty}$. Given any $c>0$, the family $X^{z}$, $\Re z>c$, is a holomorphic family in $\sL_{c^{-1}p}$.
\end{enumerate}
\end{lemma}
\begin{proof}
 The above result is standard. We include a proof for the reader's convenience. We may assume $X\neq 0$, since otherwise the result is trivially true. 
  Let $(\xi_j)_{j\geq 0}$ be an orthonormal eigenfamily of $X$ such that $X\xi_j=\lambda_j(X)\xi_j$ for all $j\geq 0$, where 
 $\lambda_0(X)\geq \lambda_1(X)\geq \cdots $ are the eigenvalues of $X$. It follows from~(\ref{eq:Com.Xz-definition}) that  
\begin{equation*}
 X^{z} = \sum_{\lambda_j(X)>0} \lambda_j(X)^{z} |\xi_j\rangle \langle\xi_j |, \qquad \Re z>0,
\end{equation*}
 where the series converges in $\sK$ uniformly on each halfspace $\Re z>c>0$. This ensures that $X^{z}$, $\Re z>0$, is a holomorphic family of compact operators. 
  
Suppose now that $X\in \sL_{p,\infty}$. Let $c>0$ and $\epsilon \in (0,c)$. Then $X^{z}=X^{c+\epsilon}X^{z-c-\epsilon}$, $\Re z>c+\epsilon$, is a holomorphic family in $\sL_{c^{-1}p}$, since $X^{c+\epsilon}$ is an operator in $\sL_{(c+\epsilon)^{-1}p,\infty}\subseteq \sL_{c^{-1}p}$ and $X^{z-c-\epsilon}$, $\Re z>c+\epsilon$, is a holomorphic family in $\sK$ by the first part. As this is true for all $\epsilon \in (0,c)$, it follows that $X^{z}$, $\Re z>c$, is a holomorphic family in $\sL_{c^{-1}p}$. 
The proof is complete.  
\end{proof}

The following lemma is a version of Lemma~\ref{lem:Tauberian.powers} for imaginary powers.  

\begin{lemma}\label{lem:Com.powers}
 If $X$ is a (non-zero) positive compact operator on $\sH$, then $X^{it}$, $t\in \R$, is a bounded strongly continuous family in $\sL(\sH)$. 
\end{lemma}
\begin{proof}
 As with Lemma~\ref{lem:Tauberian.powers}, the above result is standard. We include a proof for the reader's convenience. 
  Let $(\xi_j)_{j\geq 0}$ be an orthonormal eigenfamily of $X$ such that $X\xi_j=\lambda_j(X)\xi_j$ for all $j\geq 0$, where 
 $\lambda_0(X)\geq \lambda_1(X)\geq \cdots $ are the eigenvalues of $X$. The family $X^{it}$, $t\in \R$, then is a strongly continuous family in $\sL(\sH)$, since,  for every $\xi\in \sH$, we have
 \begin{equation*}
 X^{it}\xi= \sum_{\lambda_j(X) >0} \lambda_j(X)^{it} \scal{\xi_j}{\xi}\xi_j, \qquad t\in \R, 
\end{equation*}
where the series converges in $\sH$ uniformly on compact sets of $\R$. Note also that $\|X^{it}\|=1$, since $|\lambda_j(X)^{it}|=1$ if $\lambda_j(X)>0$ (such an eigenvalue exists, since $X\neq 0$). In particular, the family $X^{it}$, $t\in \R$, is bounded in $\sL(\sH)$.  The result is proved. 
\end{proof}

\subsection{Spectral triples}\label{subsec:Prelim.ST}
We shall use the following definition of spectral triples.

\begin{definition}
 A \emph{spectral triple} is a triple $(\sA, \sH, D)$, where $\sA$ is a unital $*$-algebra represented by bounded operators on the (separable) Hilbert space $\sH$ and $D$ is a selfadjoint unbounded operator on $\sH$ with compact resolvent such that
\begin{equation*}
a\left(\dom(D)\right)\subseteq \dom(D) \quad \text{and} \quad [D,a]\in \sL(\sH) \qquad \forall a \in \sA.
\end{equation*}
We further say that $(\sA,\sH, D)$ is \emph{$p$-summable}, with $p>0$, if the partial inverse $D^{-1}$ is in $\sL_{p,\infty}$.
\end{definition}

In what follows, we let $(\sA,\sH,D)$ be a spectral triple. We denote by $\overline{\sA}$ the closure of $\sA$ in $\sL(\sH)$, and we let $\overline{\sA}_{+}$ be its cone of positive elements. We also denote by $\overline{\sA}_{++}$ the open cone of its positive invertible elements (i.e., selfadjoint elements with positive spectrum), and we set $\sA_{++}:=\sA \cap \overline{\sA}_{++}$.

We denote by $\Pi_0$ the orthogonal projection onto $\ker D$. This is a finite-rank operator, since $D$ has compact resolvent. We define the powers $|D|^z$, $z\in \C$, as in Section~\ref{subsec:Prelim.powers}, i.e., with the convention that $|D|^z=0$ on $\ker D$. In particular, $|D|^0=1-\Pi_0$, and $|D|^{-1}$ is the partial inverse of $|D|$. We also let $F=D|D|^{-1}$ be the sign of $D$. If $(\sA,\sH,D)$ is $p$-summable, then $|D|^{-1}=|D^{-1}|$ is contained in $\sL_{p,\infty}$, and so we have
\begin{equation*}
 |D|^{-q}\in \sL_{q^{-1}p,\infty} \qquad \text{for all $q>0$}.
\end{equation*}

If $a\in \overline{\sA}_{++}$, then the action of $a$ on $\sH$ is an isomorphism. Thus, $aD^2a$ is a selfadjoint densely defined operator on $\sH$ with non-negative spectrum whose domain is $a^{-1}(\dom(D^2))$. Moreover, the operator $a(D^2+1)a$ has the compact inverse $a^{-1}(D^2+1)^{-1}a^{-1}$ and differs from $aD^2a$ by the bounded operator $a^2$, and so $aD^2a$ has compact resolvent. Thus, it has pure discrete spectrum, and we can list its \emph{positive} eigenvalues as a sequence, 
\begin{equation}\label{eq:ST.eigenvalues-aD2a}
0< \lambda_0\big(aD^2a\big) \leq  \lambda_1\big(aD^2a\big) \leq  \lambda_2\big(aD^2a\big) \leq \cdots, 
\end{equation}
where each eigenvalue is repeated according to multiplicity.

In what follows by a \emph{Fr\'echet subalgebra} of $\overline{\sA}$ we shall mean a subalgebra $\sB$ which is endowed with a Fr\'echet algebra topology which is stronger than the topology of $\overline{\sA}$. Recall that if $\sB$ is a unital Fr\'echet algebra, then the inverse map $x\rightarrow x^{-1}$ is continuous on its set of invertible elements $\sB^{-1}$ (see~\cite{Pf:PAMS, Ze:BAPS}). Therefore, if $\sB$ is a Fr\'echet subalgebra of $\overline{\sA}$ containing the unit, then it is holomorphically closed (i.e., it is closed under holomorphic functional calculus) if and only if it is inverse closed (i.e., $\sB^{-1}=\overline{\sA}^{-1}\cap \sB$).

Considering Fr\'echet subalgebras ensures that Banach--Steinhaus theorem holds. The need for holomorphic closedness stems from the following lemma.

In what follows, given a subalgebra $\sB\subseteq \overline{\sA}$, we set $\sB_{++}:=\sB\cap \overline{\sA}_{++}$. 

\begin{lemma}\label{lem:Tauberian.hol-closed-powers}
 Let $\sB$ be a holomorphically closed Fr\'echet subalgebra of $\overline{\sA}$ containing $\sA$. For every $a\in \sB_{++}$, the family $a^s$, $s\in \C$, is a holomorphic family in $\sB$. 
 \end{lemma}
\begin{proof}
Let $a\in \sB_{++}$. If $\sB=\overline{\sA}$, then the result is an immediate consequence of the power series expansion, 
 \begin{equation*}
 a^s=e^{s\log a}=\sum_{k\geq 0} \frac{1}{k!} s^k(\log a)^k,
\end{equation*}
since the above series converges normally in $\sB$ uniformly on compact subsets of $\C$. 

In general, we proceed as follows. As $a$ has positive spectrum, its spectrum is contained in some interval $[c_1,c_2]$ with $0<c_1<c_2$. Let $\Gamma$ be a direct oriented piecewise $C^1$ contour in the halfplane $\Re \lambda>0$ whose interior contains $[c_1,c_2]$. We then have
\begin{equation*}
 a^s= \frac{1}{2i\pi} \int_\Gamma \lambda^s (\lambda -a)^{-1}d\lambda, \qquad s\in \C. 
\end{equation*}
 Note that $\lambda \rightarrow (\lambda -a)^{-1}$ is a continuous map from $\Gamma$ to $\sB$. This allows us to regard the above integral as a $\sB$-valued integral.  
 
 For $\lambda \in \Gamma$ and $s\in \C$, we have
 \begin{equation*}
 \lambda^s-1-s\log \lambda = \sum_{k\geq 2} \frac{1}{k!} s^k(\log \lambda)^k = s^2h(\lambda,s), 
\end{equation*}
 where we have set 
 \begin{equation*}
 h(\lambda,s)= \sum_{k\geq 0} \frac{1}{(k+2)!} s^{k}(\log \lambda)^{k+2}, \qquad (\lambda,s)\in \Gamma\times \C. 
\end{equation*}
 Note that $h(\lambda,s)$ is a continuous function on $\Gamma\times \C$. We then have 
 \begin{align*}
 a^s-1-s\log a &=  \frac{1}{2i\pi} \int_\Gamma ( \lambda^s-1-s\log \lambda )(\lambda -a)^{-1}d\lambda\\ 
 &= \frac{s^2}{2i\pi } \int_\Gamma h(\lambda,s)(\lambda -a)^{-1}d\lambda.  
\end{align*}
 The integrand $ h(\lambda,s)(\lambda -a)^{-1}$ is a continuous map from $\Gamma\times \C$ to $\sB$, and hence it is bounded in $\sB$ on $\Gamma\times \overline{D(0,1)}$, where $\overline{D(0,1)}$ is the closed unit disk. Therefore, given any continuous semi-norm $\fp$ on $\sB$, there is a constant $C_{\fp\Gamma}>0$ such that
 \begin{equation*}
 \fp\left(a^s-1-s\log a\right)\leq \frac{|s|^2}{2\pi} \int_\Gamma \fp\left(h(\lambda,s)(\lambda -a)^{-1}\right)|d\lambda|\leq C_{\fp\Gamma}|s|^2. 
\end{equation*}
 Thus, near $s=0$, we have 
  \begin{equation*}
 a^s-1-s\log a=\op{O}(|s|^2) \qquad \text{in}\ \sB. 
\end{equation*}
This shows that the $\sB$-valued family $a^s$, $s\in \C$, is holomorphic at $s=0$ (with respect to the $\sB$-topology). The 1-parameter group property $a^s=a^{s_0}a^{s-s_0}$ further ensures it is holomorphic at every $s_0\in \C$. This proves the result. 
\end{proof}

We shall also need the following result on the commutators $[D,a^z]$, $a\in \sA_{++}$, $z\in \C$.
 
\begin{lemma}\label{lem:CSZ.Daz}
 If $a\in {\sA}_{++}$, then $[D,a^z]$, $z\in \C$, is a holomorphic family in $\sL(\sH)$. Moreover, for any closed vertical strip $\Sigma=\{c_1\leq \Re z \leq c_2\}\subseteq \C$, there is a constant $C_\Sigma(a)>0$, such that
 \begin{equation}\label{eq:App.boundedness-comm-vertical}
 \left\| [D,a^z]\right\| \leq C_{\Sigma}(a) \left(1+|\Im z|\right) \qquad \forall z\in \Sigma. 
\end{equation}
\end{lemma}
 \begin{proof}
Let $\sA_1$ be the $*$-subalgebra of $\overline{\sA}$ provided by Proposition~\ref{prop:App.hol-extension}. By Remark~\ref{rmk:App.hol-extension} it is a holomorphically closed Fr\'echet subalgebra of $\overline{\sA}$ containing $\sA$ with respect to the norm $\|T\|_1:=\|T\|+\|[D,T]\|$. Given $a\in {\sA}_{++}$, Lemma~\ref{lem:Tauberian.hol-closed-powers} then ensures that $a^z$, $z\in \C$, is a holomorphic family in $\sA_1$. As $T\rightarrow [D,T]$ is a continuous linear map from $\sA_1$ to $\sL(\sH)$, it follows that $[D,a^z]$, $z\in \C$, is a holomorphic family in $\sL(\sH)$. 

The above property implies that, for any compact $K\subseteq \C$, the family $[D,a^z]$, $z\in K$, is bounded in $\sL(\sH)$. Moreover, if $z=x+iy$ with $x,y\in \R$, then 
\begin{equation*}
  \big[D,a^z] =  \big[D,a^x]a^{iy} +  a^x\big[D,a^{iy}]. 
\end{equation*}
We observe that if $x$ remains in a closed interval $[\alpha,\beta]$ 
and $y$ varies over $\R$, then $a^x$ and $a^{iy}$ remain in bounded subsets of $\sL(\sH)$ and $[D,a^x]$ remains in a bounded subset of $\sL(\sH)$. Therefore, in order to prove the estimates~(\ref{eq:App.boundedness-comm-vertical}) for closed vertical strips it is enough to do it for the imaginary line $\Sigma =i\R$. 

As $i[-1,1]$ is a compact subset of $\C$, there is $C>0$ such that
\begin{equation*}
 \left\| [D,a^{it}]\right\| \leq C \qquad \forall t\in [-1,1]. 
\end{equation*}
If $y\geq 1$ and we write $y=t+\ell$ with $t\in [0,1]$ and $\ell \in \N$, then 
\begin{align*}
   [D,a^{iy}]  & = a^{it}   [D,a^{i\ell} ] +   [D,a^{it}] a^{i\ell} \\
   & = \sum_{0\leq j \leq \ell-1} a^{i(t+j)}   [D,a^{i} ]a^{i(\ell-1-j)} +  [D,a^{it}] a^{i\ell}. 
\end{align*}
Thus,
\begin{equation*}
  \left\| [D,a^{iy}]\right\| \leq \ell  \left\| [D,a^{i}]\right\| +  \left\| [D,a^{it}]\right\|\leq C(\ell+1)\leq C(|y|+1). 
\end{equation*}
Likewise, if $y\leq -1$, then  $ \| [D,a^{iy}]\| \leq C(1+|y|)$. This proves~(\ref{eq:App.boundedness-comm-vertical}) for $\Sigma=i\R$. The proof is complete. 
\end{proof}
 
\section{Double Operator Integrals}\label{chap:DOIs}
This section gathers the results on double operator integrals (DOIs) that are needed for the proof of the Factorization Lemma (Theorem~\ref{thm:Intro.factorization}) in Section~\ref{sec:Factorization}. These results rely on basic techniques of DOIs that were established by Birman and Solomyak in the 60s and 70s. Throughout this section we let $(\sA,\sH,D)$ be a spectral triple. 

The proof of Theorem~\ref{thm:Intro.factorization} relies on expressing commutators in terms of operators of the form, 
\begin{equation}\label{eq:Com.Psiadd'}
 \Psi_{\delta,\delta'}^\alpha(T) = \int_0^\infty \mu^{1+\alpha}|D|^\delta\left(D^2+\mu^2\right)^{-1}T |D|^{\delta'}\left(D^2+\mu^2\right)^{-1}d\mu, \qquad T\in \sL(\sH),
\end{equation}
where $\alpha \in (-2,2)$ and $\delta,\delta'\in [0,2]$ satisfy one of the following two conditions:
\begin{enumerate}
 \item[(i)] $\alpha+\delta+\delta'<2$, or
 
 \item[(ii)] $\alpha+\delta+\delta'=2$ and $\delta,\delta'\in (0,2)$. 
\end{enumerate}
These operators (when they make sense) are special instances of double operator integrals (DOIs) in the sense of Birman--Solomyak~\cite{BS:PMF66, BS:PMF73,BS:IEOT03} (see also~\cite{ST:LNM19}). In fact, in the case (ii), up to change of variable and scaling, the operator $ \Psi_{\delta,\delta'}^\alpha(T) $ is the operator $\Theta^{A_0,A_1}_{\frac{1}{2}\delta,\frac{1}{2}\delta'}(T)$ with $A_0=A_1=D^2$ from~\cite[Section 1]{BS:JSM92}. Therefore, at least in the case  $\alpha+\delta+\delta'=2$, the properties of our operators $\Psi_{\delta,\delta'}^\alpha(T)$ follow from standard properties of DOIs as described in~\cite{BS:PMF66, BS:PMF73,BS:IEOT03, ST:LNM19}. 
 
 For our purpose we only need the boundedness of DOIs on $\sL(\sH)$.  Moreover, it is important to incorporate condition~(i) as well, even though it is simpler. For these reasons, we provide a unified account that focuses on boundedness on $\sL(\sH)$ under our two conditions, and includes the original arguments of~\cite{BS:PMF73,BS:IEOT03} for our class of DOIs under condition~(ii).  
 
 Actually, the main issue at stake is to make sense of the integrals appearing in~(\ref{eq:Com.Psiadd'}). To ease the exposition we introduce some notation. For $0\leq \delta \leq 2$, set 
 \begin{equation*}
 R_\delta(\mu)=|D|^\delta\left(D^2+\mu^2\right)^{-1}, \qquad \mu\geq 0.  
\end{equation*}
with the convention that $|D|^0=1-\Pi_0$. Note that $R_\delta(\mu)$ is a bounded operator which is compact for $0\leq \delta<2$. For $0\leq \delta,\delta' \leq 2$, we also define
 \begin{align}
 \Delta_{\delta,\delta'}(T;\mu):= &  R_\delta(\mu)TR_{\delta'}(\mu)\nonumber \\
=& |D|^\delta\left(D^2+\mu^2\right)^{-1}T |D|^{\delta'}\left(D^2+\mu^2\right)^{-1}, \quad \mu\geq 0, \ T\in \sL(\sH). 
\label{eq:Com.Deltadd'mu} 
\end{align}
Note that $\Delta_{\delta,\delta'}(T;\mu)$ is always a bounded operator and is compact whenever $\min(\delta,\delta')<2$.

 \begin{proposition} \label{lem:Com.estimate-Delta}
 Let $\alpha \in (-2,2)$ and $\delta,\delta'\in [0,2]$ be such that $\alpha+\delta+\delta'<2$. Then, for every $T\in \sL(\sH)$, the integral~(\ref{eq:Com.Psiadd'}) converges in $\sK$ in Bochner's sense. This defines a bounded operator $T\rightarrow\Psi^\alpha_{\delta,\delta'}(T)$ from $\sL(\sH)$ to $\sK$.
 \end{proposition}
 \begin{proof} 
We first show that there is $C_{\alpha \delta \delta'}>0$ such that
 \begin{equation}\label{eq:Com.estimate-Delta1}
 \int_0^\infty\mu^{1+\alpha} \left\| \Delta_{\delta,\delta'}(T;\mu)\right\| d\mu \leq C_{\alpha \delta \delta'}\|T\| \qquad \forall T\in \sL(\sH). 
\end{equation}
Let $T\in \sL(\sH)$. We have 
\begin{equation}\label{eq:Com.estimate-int-Delta1}
  \int_0^\infty  \mu^{1+\alpha} \left\| \Delta_{\delta,\delta'}(T;\mu)\right\| d\mu \leq \|T\|  \int_0^\infty  \mu^{1+\alpha} \| R_{\delta}(\mu)\| \|R_{\delta'}(\mu)\|d\mu. 
\end{equation}
 Moreover, given any $\mu>0$ and $\delta \in [0,2]$, we have
\begin{equation*}
\left \|R_\delta(\mu)\right\|= \mu^{\delta-2}\big\||(\mu^{-1}D)|^{\delta}\big((\mu^{-1}D)^2+1\big)^{-1}\big\|\leq \mu^{-2+\delta}\max_{t\geq 0} t^{\delta}(t^2+1)^{-1}.
\end{equation*}
In particular, we see that
\begin{equation}\label{eq:Com.Rd-bigO}
 \left  \|R_\delta(\mu)\right\|=\op{O}\left(\mu^{-2+\delta}\right)\qquad \text{as}\ \mu \rightarrow \infty. 
\end{equation}
In addition, as $|D|^\delta$ vanishes on $\ker D$, we have $\|R_\delta(\mu)\|\leq \||D|^{\delta-2}\|$ for all $\mu\geq 0$. In particular, $\|R_\delta(\mu)\|$ remains bounded as $\mu\rightarrow 0^+$. 
 
 As $\alpha+\delta+\delta'<2$, the estimate~(\ref{eq:Com.Rd-bigO}) ensures that $\mu^{1+\alpha}\| R_{\delta}(\mu)\| \|R_{\delta'}(\mu)\|=\op{O}(\mu^r)$, where 
 \begin{equation*}
 r=(1+\alpha)+(-2+\delta)+(-2+\delta')=-1-[2-(\alpha+\delta+ \delta')]<-1.
\end{equation*}
 Together with the boundedness of $\|R_\delta(\mu)\|$ and $\|R_{\delta'}(\mu)\|$ near $\mu=0$ and the fact that $\alpha>-2$, i.e.,  $1+\alpha>-1$, we deduce that
 \begin{equation*}
  \int_0^\infty  \mu^{1+\alpha} \| R_{\delta}(\mu)\| \|R_{\delta'}(\mu)\|d\mu <\infty.
\end{equation*}
Combining this with~(\ref{eq:Com.estimate-int-Delta1}) gives the estimate~(\ref{eq:Com.estimate-Delta1}). 

The estimate~(\ref{eq:Com.estimate-Delta1}) implies that the map $\mu \rightarrow \mu^{1+\alpha}\Delta_{\delta,\delta'}(T;\mu)$ is integrable in $\sL(\sH)$ in Bochner's sense, and so the integral~(\ref{eq:Com.Psiadd'}) converges in $\sL(\sH)$ and $\|\Psi^\alpha_{\delta,\delta'}(T)\|\leq C_{\alpha\delta\delta'}\|T\|$. Moreover, as $\alpha>-2$, we have $\delta+\delta'<2-\alpha<4$, and hence $\min(\delta,\delta')<2$. Thus, one of the operators $R_\delta(\mu)$ or $R_{\delta'}(\mu)$ is compact, and so $\Delta_{\delta,\delta'}(T;\mu)$ is compact for all $\mu\geq 0$. As $\sK$ is closed in $\sL(\sH)$, the integral actually converges in $\sK$. The proof is complete. 
\end{proof}

\begin{remark}
 The arguments of the proof of the estimate~(\ref{eq:Com.estimate-Delta1}) are reminiscent of standard arguments that are used to show that spectral triples give rise to bounded Fredholm modules (see, e.g., \cite{BJ:CRAS83, CP:CJM98, CM:CRAS86, SWW:CRAS98}).  
 \end{remark}

 \begin{proposition}[Birman--Solomyak~\cite{BS:PMF66,BS:PMF73,BS:IEOT03}]  \label{lem:Com.estimate-Delta-weak}
 Let $\alpha \in (-2,2)$ and $\delta,\delta'\in (0,2)$ be such that $\alpha+\delta+\delta'=2$. Then the integral~(\ref{eq:Com.Psiadd'}) converges in the weak operator topology, i.e., for every $T\in\sL(\sH)$ there is a unique $\Psi^\alpha_{\delta,\delta'}(T)\in\sL(\sH)$ such that
 \begin{equation}\label{eq:Com.def-Psidd'-weak}
 \scal{\Psi_{\delta,\delta'}^\alpha(T)\xi}{\eta}=\int_0^\infty\mu^{1+\alpha} \bigscal{\Delta_{\delta,\delta'}(T;\mu)\xi}{\eta} d\mu \qquad \forall \xi,\eta\in \sH. 
 \end{equation}
 This defines a bounded operator $T\rightarrow\Psi^\alpha_{\delta,\delta'}(T)$ from $\sL(\sH)$ to itself.
 \end{proposition}

\begin{remark}
The above result is a special case of a general result for a larger class of DOIs for which we have separation of variables (see~\cite{BS:PMF73, BS:IEOT03, Pe:FAA85, ST:LNM19}). For the reader's convenience, a self-contained proof of Proposition~\ref{lem:Com.estimate-Delta-weak} is given in Appendix~\ref{app:BS}.  
\end{remark}

We will need the following factorization properties. 

 \begin{lemma}\label{lem:Com.factorization-Psidd'}
Suppose that $\alpha \in (-2,2)$ and  $\delta,\delta',\epsilon,\epsilon'\in [0,2]$ satisfy one of the following:
\begin{itemize}
 \item[(a)] $\alpha+\delta+\delta'+\epsilon+\epsilon'<2$ with $\max(\delta+\epsilon, \delta'+\epsilon')\leq 2$.
 \item[(b)]  $\alpha+\delta+\delta'+\epsilon+\epsilon'=2$ with $\delta+\epsilon$ and  $\delta'+\epsilon'$ in $(0,2)$.
\end{itemize}
Then, we have
 \begin{equation}\label{eq:Com.factorization-Psidd'}
 \Psi_{\delta,\delta'}^\alpha(T) = |D|^{-\epsilon} \Psi_{\delta+\epsilon,\delta'+\epsilon'}^\alpha(T) |D|^{-\epsilon'} \qquad \forall T\in \sL(\sH). 
\end{equation}
\end{lemma}
\begin{proof}
 Let $T\in \sL(\sH)$. Note that $R_\delta(\mu)=|D|^{-\epsilon}R_{\delta+\epsilon}(\mu)$. Likewise, $R_{\delta'}(\mu)=|D|^{-\epsilon'}R_{\delta'+\epsilon'}(\mu)=R_{\delta'+\epsilon'}(\mu)|D|^{-\epsilon'}$. Thus, 
\begin{equation}\label{eq:Com.factorization-Deltadd'}
 \Delta_{\delta,\delta'}(T;\mu)=  |D|^{-\epsilon}  \Delta_{\delta+\epsilon,\delta'+\epsilon'}(T;\mu)|D|^{-\epsilon'}, \qquad \mu\geq 0. 
\end{equation}

If $\alpha +\delta+\delta'+\epsilon+\epsilon'<2$ and $\max(\delta+\epsilon, \delta'+\epsilon')\leq 2$, then  $\Psi_{\delta,\delta'}^\alpha(T)$ and $\Psi_{\delta+\epsilon,\delta'+\epsilon'}^\alpha(T)$ are both given by Bochner integrals, and so we have
\begin{align*}
 \Psi_{\delta,\delta'}^\alpha(T) & = \int_0^\infty \mu^{1+\alpha} |D|^{-\epsilon}  \Delta_{\delta+\epsilon,\delta'+\epsilon'}(T;\mu)|D|^{-\epsilon'}d\mu \\
 &=|D|^{-\epsilon} \left( \int_0^\infty \mu^{1+\alpha}  \Delta_{\delta+\epsilon,\delta'+\epsilon'}(T;\mu)d\mu\right)|D|^{-\epsilon'}\\
 & = |D|^{-\epsilon} \Psi_{\delta+\epsilon,\delta'+\epsilon'}^\alpha(T) |D|^{-\epsilon'}. 
\end{align*}

Suppose now that $\alpha+\delta+\delta'+\epsilon+\epsilon'=2$ with $\delta+\epsilon$ and  $\delta'+\epsilon'$ in $(0,2)$. Let $\xi ,\eta\in \sH$. In view of~(\ref{eq:Com.def-Psidd'-weak}) we have
\begin{align*}
 \bigscal{ |D|^{-\epsilon} \Psi_{\delta+\epsilon,\delta'+\epsilon'}^\alpha(T) |D|^{-\epsilon'}\xi}{\eta} & = 
  \bigscal{ \Psi_{\delta+\epsilon,\delta'+\epsilon'}^\alpha(T)|D|^{-\epsilon'}\xi}{|D|^{-\epsilon} \eta}\\
  &= \int_0^\infty\mu^{1+\alpha} \bigscal{\Delta_{\delta+\epsilon,\delta'+\epsilon'}(T;\mu)|D|^{-\epsilon'}\xi}{|D|^{-\epsilon}\eta} d\mu. 
\end{align*}
Thanks to~(\ref{eq:Com.factorization-Deltadd'}) we have
\begin{equation*}
  \bigscal{\Delta_{\delta+\epsilon,\delta'+\epsilon'}(T;\mu)|D|^{-\epsilon'}\xi}{|D|^{-\epsilon}\eta} =
   \bigscal{|D|^{-\epsilon}\Delta_{\delta+\epsilon,\delta'+\epsilon'}(T;\mu)|D|^{-\epsilon'}\xi}{\eta}=  
   \bigscal{\Delta_{\delta,\delta'}(T;\mu)\xi}{\eta}. 
\end{equation*}
Thus, 
\begin{equation*}
  \bigscal{ |D|^{-\epsilon} \Psi_{\delta+\epsilon,\delta'+\epsilon'}^\alpha(T) |D|^{-\epsilon'}\xi}{\eta} =  
  \int_0^\infty\mu^{1+\alpha} \bigscal{\Delta_{\delta,\delta'}(T;\mu)\xi}{\eta} d\mu = 
   \bigscal{  \Psi_{\delta,\delta'}^\alpha(T)\xi}{\eta}. 
\end{equation*}
This gives the equality~(\ref{eq:Com.factorization-Psidd'}). The proof is complete. 
\end{proof}

\begin{remark}
 As $|D|^0=1-\Pi_0$, the factorization~(\ref{eq:Com.factorization-Psidd'}) for $\epsilon=\epsilon'=0$ means that, for all $T\in \sL(\sH)$, we have
 \begin{equation*}
 (1-\Pi_0)\Psi^\alpha_{\delta,\delta'}(T)(1-\Pi_0)=\Psi^\alpha_{\delta,\delta'}(T).
\end{equation*}
That is, $\ker D\subseteq \ker (\Psi^\alpha_{\delta,\delta'}(T))$ and $\ran(\Psi^\alpha_{\delta,\delta'}(T))\subseteq (\ker D)^\perp=\overline{\ran D}$. 
\end{remark}

 Assume further that $(\sA, \sH, D)$ is $p$-summable, so that $|D|^{-q}\in \sL_{q^{-1}p,\infty}$ for all $q>0$. We then have the following consequence of Lemma~\ref{lem:Com.factorization-Psidd'}. 
 
\begin{lemma}\label{lem:Com.Psidd'.sLr}
 Suppose that $\alpha+\delta+\delta'<2$. 
 \begin{enumerate}
 \item[(i)] If $\max(\delta,\delta')<2$, then $\Psi^{\alpha}_{\delta,\delta'}$ gives rise to a bounded operator, 
\begin{equation*}
 \Psi_{\delta,\delta'}^\alpha: \sL(\sH)\longrightarrow \sL_{r,\infty}, \qquad \textup{with}\ r=(2-\alpha-\delta-\delta')^{-1}p. 
\end{equation*}

\item[(ii)] If $\max(\delta,\delta')=2$, then the same result holds for any $r>(2-\alpha-\delta-\delta')^{-1}p$.   
\end{enumerate}
\end{lemma}
\begin{proof}
 Suppose that $\max(\delta,\delta')<2$. The function $(\epsilon,\epsilon')\rightarrow \alpha+\delta+\delta'+\epsilon+\epsilon'$ is continuous and maps $(0,2-\delta)\times (0,2-\delta')$ onto $(\alpha+\delta+\delta',\alpha+4)$. As $\alpha+\delta+\delta'<2$ and $\alpha>-2$, the latter interval contains $2$, and so we can always find $(\epsilon,\epsilon')\in (0,2-\delta)\times (0,2-\delta')$ with $\alpha+\delta+\delta'+\epsilon+\epsilon'=2$. Then $\delta+\epsilon$ and $\delta'+\epsilon'$ are in $(0,2)$. By Lemma~\ref{lem:Com.factorization-Psidd'} we then have
 \begin{equation*}
 \Psi_{\delta,\delta'}^\alpha(T) = |D|^{-\epsilon} \Psi_{\delta+\epsilon,\delta'+\epsilon'}^\alpha(T) |D|^{-\epsilon'},  \qquad T\in \sL(\sH). 
\end{equation*}
Here $ \Psi_{\delta+\epsilon,\delta'+\epsilon'}^\alpha:\sL(\sH)\rightarrow \sL(\sH)$ is bounded since  $\max(\delta+\epsilon, \delta'+\epsilon')<2$. Moreover, the operators $|D|^{-\epsilon}$ and $|D|^{-\epsilon'}$ are contained in $\sL_{\epsilon^{-1}p,\infty}$ and  $\sL_{(\epsilon')^{-1}p,\infty}$, respectively. By using H\"older's inequality (Proposition~\ref{prop:Hoelder}) we then deduce that $\Psi_{\delta,\delta'}^\alpha$ maps continuously $\sL(\sH)$ to $\sL_{r,\infty}$, with  
$r=(\epsilon+\epsilon')^{-1}p=(2-\alpha-\delta-\delta')^{-1}p$. This proves the first part. 

Suppose that $\max(\delta, \delta')=2$. Assume that $\delta'=2$. The inequality $\alpha+\delta+\delta'<2$ then means that $\alpha+\delta<0$. Let $\epsilon \in (0,-(\alpha+\delta)]$, and set $\epsilon_1=-(\alpha+\delta)-\epsilon$. Then $\epsilon_1\in [0,-(\alpha+\delta))\subseteq [0,2)$ and $\alpha+\delta+\delta'+\epsilon_1=2-\epsilon<2$. Therefore, by Lemma~\ref{lem:Com.factorization-Psidd'} we have 
\begin{equation*}
  \Psi_{\delta,\delta'}^\alpha(T) = |D|^{-\epsilon_1} \Psi_{\delta+\epsilon_1,\delta'}^\alpha(T) ,  \qquad T\in \sL(\sH). 
\end{equation*}
As $|D|^{-\epsilon_1}\in \sL_{\epsilon_1^{-1}p,\infty}$, in the same way as above we see that   
$\Psi_{\delta,\delta'}^\alpha$ maps continuously $\sL(\sH)$ to $\sL_{r,\infty}$, with  
$r=\epsilon_1^{-1}p= (-\alpha-\delta-\epsilon)^{-1}p=(2-\alpha-\delta-\delta'-\epsilon)^{-1}p$. Here $\epsilon$ ranges over $(0,-(\alpha+\delta)]$, and so $r$ ranges 
over $((2-\alpha-\delta-\delta')^{-1}p,\infty)$. The same result holds in case $\delta=2$. This gives the 2nd part and completes the proof.  
\end{proof}

\begin{remark}\label{rmk:Com.Psidd'.sLr2}
 Suppose that $\alpha+\delta+\delta'=2$ and  $\delta,\delta' \in (0,2)$. A well-known result of the theory of DOI (see~\cite{BS:PMF73,BS:IEOT03}) asserts that DOIs like $\Psi_{\delta,\delta'}^\alpha$ induce bounded operators on $\sL_1$. By interpolation, we then obtain the boundedness of $\Psi_{\delta,\delta'}^\alpha$ on all weak Schatten classes $\sL_{q,\infty}$ with $q>1$ (see~\cite{BS:JSM92, BS:IEOT03,ST:LNM19}). 
\end{remark}

\begin{remark}
 Suppose that $\alpha+\delta+\delta'<2$. By using Remark~\ref{rmk:Com.Psidd'.sLr2} and arguing as in the proof of Lemma~\ref{lem:Com.Psidd'.sLr} we obtain that if $\max(\delta, \delta')<2$, then, given any $q>1$, the operator $\Psi_{\delta,\delta'}^{\alpha}$ induces a bounded operator, 
\begin{equation*}
 \Psi_{\delta,\delta'}^\alpha: \sL_{q,\infty}\longrightarrow \sL_{r,\infty}, \qquad r^{-1}=q^{-1}+(2-\alpha-\delta-\delta')p^{-1}. 
\end{equation*}
If $\max(\delta,\delta')=2$, the same result holds for any $r>(q^{-1}+(2-\alpha-\delta-\delta')p^{-1})^{-1}$.
\end{remark}

\section{Factorization Lemma}\label{sec:Factorization}
This section is devoted to proving the Factorization Lemma (Theorem~\ref{thm:Intro.factorization}), which provides factorization formulas for fractional commutators. These formulas play a crucial role in the proof of the Enhanced CSZ Lemma (Theorem~\ref{thm:Intro.CSZ}) in Section~\ref{chap:CSZ}, and hence in the proof of the Tauberian Theorem (Theorem~\ref{thm:Intro.Tauberian}) in Section~\ref{chap:Tauberian}. In Part~I~\cite{Po:Part1} they are also used to derive commutator estimates.

Throughout this section we let $(\sA,\sH,D)$ be a $p$-summable spectral triple. We also set
\begin{equation*}
 c(\alpha) = \frac{2}{\pi} \sin \left(\frac{\pi \alpha}{2}\right), \qquad \alpha \geq 0. 
\end{equation*}

\begin{lemma} \label{lem:Com.Pi0Com|D|m-Psi}
Let $\alpha \in [0,1]$. Given any $a\in \sA$, we have 
 \begin{equation}\label{eq:Com.Pi0-Com-|D|mleq0}
  (1-\Pi_0)\left[ |D|^{-\alpha},a\right](1-\Pi_0)  =-c(\alpha)F\Psi_{1,0}^{-\alpha}\left([D,a]\right) - c(\alpha)\Psi_{0,1}^{-\alpha}\left([D,a]\right)F. 
\end{equation}
 \end{lemma}
\begin{proof}
If $\alpha=0$ then $c(\alpha)=0$, and so the r.h.s.\ of~(\ref{eq:Com.Pi0-Com-|D|mleq0}) vanishes. The l.h.s.\ vanishes as well, since, as $|D|^0=1-\Pi_0$, we have
\begin{equation*}
 (1-\Pi_0)\left[ |D|^{0},a\right](1-\Pi_0)= (1-\Pi_0)\left[ (1-\Pi_0),a\right](1-\Pi_0)=0. 
\end{equation*}
Thus, Eq.~(\ref{eq:Com.Pi0-Com-|D|mleq0}) holds trivially for $\alpha=0$. 

Suppose that $0 < \alpha \leq 1$.  If $x>0$ and $0<\beta <1$, then
\begin{equation*}
 x^{-\beta} = \frac{\sin (\pi \beta)}{\pi} \int_0^\infty \mu^{-\beta}(x+\mu)^{-1}d\mu. 
\end{equation*}
This formula is well-known. Recall it follows from the identities, 
\begin{equation*}
 \int_0^\infty \mu^{-\beta} (1+\mu)^{-1}d\mu = \int_0^1(1-\lambda)^{-\beta}\lambda^{\beta-1}d\lambda =\Gamma(1-\beta)\Gamma(\beta)= \frac{\pi}{\sin (\pi \beta)},
\end{equation*}
the last identity being Euler's reflection formula. Thus, if $x\neq 0$, and we take $\beta=\alpha/2$, then
\begin{equation*}
 |x|^{-\alpha}=(x^2)^{-\frac{\alpha}{2}}=  \frac{\sin (\pi \alpha/2)}{\pi}\int_0^\infty \mu^{-\frac{\alpha}{2}}(x^2+\mu)^{-1}d\mu 
 =c(\alpha) \int_0^\infty \mu^{1-\alpha}(x^2+\mu^2)^{-1}d\mu. 
\end{equation*}
It then follows that
\begin{equation*}
 |D|^{-\alpha} = c(\alpha) \int_0^\infty \mu^{1-\alpha}(1-\Pi_0)\left(D^2+\mu^2\right)^{-1}d\mu.  
\end{equation*}
Thus, 
\begin{align}
 (1-\Pi_0)\left[ |D|^{-\alpha},a\right](1-\Pi_0) & =  c(\alpha) \int_0^\infty \mu^{1-\alpha} (1-\Pi_0)\big[(1-\Pi_0) \left(D^2+\mu^2\right)^{-1}\!\!,a\big](1-\Pi_0)d\mu  \nonumber \\
 &= c(\alpha) \int_0^\infty \mu^{1-\alpha} (1-\Pi_0)\left[\left(D^2+\mu^2\right)^{-1},a\right](1-\Pi_0)d\mu. 
 \label{eq:Com.Com-Da-resolvent}
\end{align}

Given any $\mu>0$,  by~\cite[Lemma~2.3]{CP:CJM98} we have
\begin{align*}
\big[\left(D^2+\mu^2\right)^{-1},a\big] = &  
 -D\left(D^2+\mu^2\right)^{-1}[D,a] \left(D^2+\mu^2\right)^{-1}\\
 &  - \left(D^2+\mu^2\right)^{-1}[D,a] D\left(D^2+\mu^2\right)^{-1}. 
\end{align*}
Note that, as it is not assumed that $a(\dom (D^2))\subseteq \dom (D^2)$, a bit of care is needed in order to get the above equality (see~\cite{CP:CJM98}). Thus, in terms of the operators $\Delta_{\delta,\delta'}(T;\mu)$ defined in~(\ref{eq:Com.Deltadd'mu}) (and recalling that $|D|^0=1-\Pi_0$) we get 
\begin{align*}
 (1-\Pi_0)\left[a, \left(D^2+\mu^2\right)^{-1}\right](1-\Pi_0) = & \  F|D|\left(D^2+\mu^2\right)^{-1}[D,a] \left(D^2+\mu^2\right)^{-1}|D|^0\\
 & + |D|^0\left(D^2+\mu^2\right)^{-1}[D,a] F |D|\left(D^2+\mu^2\right)^{-1}\\
 = & \ F\Delta_{1,0}\left([D,a];\mu\right) + \Delta_{0,1}\left([D,a];\mu\right)F. 
\end{align*}
Combining this with~(\ref{eq:Com.Com-Da-resolvent}) then gives 
\begin{align}
 (1-\Pi_0)\left[ |D|^{-\alpha},a\right](1-\Pi_0)  = & \ - c(\alpha) \int_0^\infty \mu^{1-\alpha} F\Delta_{1,0}\left([D,a];\mu\right)d\mu \nonumber\\
 & -c(\alpha) \int_0^\infty \mu^{1-\alpha} \Delta_{0,1}\left([D,a];\mu\right)Fd\mu \\ 
 = &\ -c(\alpha)F\Psi_{1,0}^{-\alpha}\left([D,a]\right) - c(\alpha)\Psi_{0,1}^{-\alpha}\left([D,a]\right)F, \nonumber
\end{align}
where the integrals converge in $\sL(\sH)$ in Bochner's sense (\emph{cf}.~Proposition~\ref{lem:Com.estimate-Delta}). This proves the result for $0<\alpha\leq 1$. The proof is complete. 
\end{proof}

We are now in a position to prove Theorem~\ref{thm:Intro.factorization}.

\begin{proof}[Proof of Theorem~\ref{thm:Intro.factorization}]
Assume that the spectral triple $(\sA, \sH, D)$ is $p$-summable with $p>0$ arbitrary. Let $\alpha\in [0,1]$. We have to show that if $\beta,\gamma\in [0,1]$ are such that, either $\beta+\gamma<\alpha+1$, or $\beta+\gamma=\alpha+1$ with $\max(\beta,\gamma)<1$, then there is a bounded operator $\Phi_{\beta,\gamma}^{\alpha}:\sL(\sH)\rightarrow \sL(\sH)$ such that, for all $a\in \sA$, we have
 \begin{equation}\label{eq:Com.factorization}
  \left[ |D|^{-\alpha},a\right]  =|D|^{-\beta}\Phi_{\beta,\gamma}^{\alpha}\left([D,a]\right)|D|^{-\gamma} 
  + |D|^{-(\alpha+1)}F[D,a]\Pi_0 +\Pi_0[D,a]F|D|^{-(\alpha+1)}. 
\end{equation}

Let $a\in \sA$. As $\Pi_0|D|^{-\alpha}=|D|^{-\alpha}\Pi_0=0$, we have
 \begin{align}
 \left[ |D|^{-\alpha},a\right]&=(1-\Pi_0)\left[ |D|^{-\alpha},a\right](1-\Pi_0)+ \Pi_0 \left[ |D|^{-\alpha},a\right] + (1-\Pi_0)\left[ |D|^{-\alpha},a\right]\Pi_0 \nonumber \\
 &= (1-\Pi_0)\left[ |D|^{-\alpha},a\right](1-\Pi_0)+ |D|^{-\alpha}a\Pi_0 -\Pi_0a|D|^{-\alpha} .
  \end{align}
 Note that
 \begin{equation}
 |D|^{-\alpha}a\Pi_0=|D|^{-(\alpha+1)}FDa\Pi_0=|D|^{-(\alpha+1)}F[D,a]\Pi_0.
\end{equation}
Likewise, we have 
\begin{equation*}
 -\Pi_0a|D|^{-\alpha}= -\Pi_0aDF|D|^{-(\alpha+1)}= \Pi_0[D,a]F|D|^{-(\alpha+1)}.
\end{equation*}
Thus,
\begin{equation}\label{eq:Com.Com-|D|ml-Pi0F}
 \left[ |D|^{-\alpha},a\right]=(1-\Pi_0)\left[ |D|^{-\alpha},a\right](1-\Pi_0)+ |D|^{-(\alpha+1)}F[D,a]\Pi_0 + \Pi_0[D,a]F|D|^{-(\alpha+1)}. 
\end{equation}

By Lemma~\ref{lem:Com.Pi0Com|D|m-Psi} we have
\begin{equation*}
 (1-\Pi_0)\left[|D|^{-\alpha},a\right](1-\Pi_0) = - c(\alpha)F\Psi_{1,0}^{-\alpha}\left([D,a]\right) - c(\alpha)\Psi_{0,1}^{-\alpha}\left([D,a]\right)F. 
\end{equation*}
Let $\beta,\gamma\in [0,1]$ be such that, either $\beta+\gamma<\alpha+1$, or $\beta+\gamma=\alpha+1$ with $\max(\beta,\gamma)<1$. Then either $-\alpha+(\beta+1)+\gamma<2$, or $-\alpha+(\beta+1)+\gamma=2$ with $\max(1+\beta,\gamma)<2$. In either case by Lemma~\ref{lem:Com.factorization-Psidd'} we have
\begin{equation*}
 \Psi^{-\alpha}_{1,0}(T)=|D|^{-\beta}\Psi^{-\alpha}_{\beta+1,\gamma}(T)|D|^{-\gamma}, \qquad T\in \sL(\sH). 
\end{equation*}
Likewise, we have 
\begin{equation*}
 \Psi^{-\alpha}_{0,1}(T)=|D|^{-\beta}\Psi^{-\alpha}_{\beta,\gamma+1}(T)|D|^{-\gamma}, \qquad T\in \sL(\sH). 
\end{equation*}
Therefore, we may write 
\begin{equation*}
 (1-\Pi_0)\left[|D|^{-\alpha},a\right](1-\Pi_0) =|D|^{-\beta}\Phi^{\alpha}_{\beta,\gamma}\left([D,a]\right)|D|^{-\gamma},
\end{equation*}
where $\Phi_{\beta,\gamma}^{\alpha}:\sL(\sH)\rightarrow \sL(\sH)$ is the bounded linear operator defined by
\begin{equation}\label{eq:Com.Phimbgam1}
 \Phi_{\beta,\gamma}^{\alpha}(T)=-c(\alpha)F\Psi^{-\alpha}_{\beta+1,\gamma}(T)-c(\alpha)\Psi^{-\alpha}_{\beta,\gamma+1}(T)F,  \qquad T\in \sL(\sH). 
\end{equation}
Combining this with~(\ref{eq:Com.Com-|D|ml-Pi0F}) yields the equality~(\ref{eq:Com.factorization}). 

Finally, to get the 2nd part of Theorem~\ref{thm:Intro.factorization} suppose that $\beta+\gamma<\alpha+1$, i.e., $-\alpha+\beta+\gamma+1<2$. If $\max(\beta, \gamma)<1$, then $\max(\beta+1,\gamma)$ and $\max(\beta,\gamma+1)$ are~$<2$. Lemma~\ref{lem:Com.Psidd'.sLr} then ensures that $\Psi^{-\alpha}_{\beta+1,\gamma}$ and $\Psi^{-\alpha}_{\beta,\gamma+1}$ both map continuously $\sL(\sH)$ to $\sL_{r,\infty}$ with $r=(2+\alpha-\beta-\gamma -1)^{-1}p=(1+\alpha-\beta-\gamma)^{-1}p$, and so $\Phi_{\beta,\gamma}^{\alpha}$ is a continuous operator from $\sL(\sH)$ to $\sL_{r,\infty}$ as well. If $\max(\beta, \gamma)=1$, then the same result holds for any 
$r>(1+\alpha-\beta-\gamma)^{-1}p$. This gives the 2nd part. The proof of Theorem~\ref{thm:Intro.factorization} is complete.
\end{proof}

\begin{remark}
 Theorem~\ref{thm:Intro.factorization} is stated in much more generality than what is needed in this paper and in Part~I~\cite{Po:Part1}. In fact, for our purpose we only need the result in the following cases:
\begin{itemize}
 \item[(i)] $0<\alpha\leq 1$ and $\beta=\gamma =0$. 
 
 \item[(ii)] $0 <\alpha\leq 1$, $\gamma =\alpha$ and $0<\beta<1$.  
\end{itemize}
The first case is used in Part~I to derive commutator estimates. The second case is needed in the proof of the Enhanced CSZ Lemma (Theorem~\ref{thm:Intro.CSZ}), for which it is absolutely crucial to be able to take $\gamma=\alpha$, including for $\alpha=1$ (see Eqs.~(\ref{eq:CSZ.Fst1}) and~(\ref{eq:CSZ.Fst2}) and Remark~\ref{rmk:CSZ.alpha} on this point). 
\end{remark}

\section{Enhanced CSZ Lemma}\label{chap:CSZ}
This section is devoted to the proofs of the Enhanced CSZ Lemma (Theorem~\ref{thm:Intro.CSZ}), of its holomorphic version (Theorem~\ref{thm:Intro.CSZhol}), and of Corollary~\ref{cor:Intro.zeta}. The latter is a key ingredient in the proof of the Tauberian Theorem (Theorem~\ref{thm:Intro.Tauberian}) in Section~\ref{chap:Tauberian}.

Throughout this section we let $(\sA, \sH, D)$ be a $p$-summable spectral triple with $p>0$ arbitrary.

\subsection{Auxiliary lemmas}\label{subsec:CSZ.aux}
In this subsection, we gather several lemmas that will be used in the proof of Theorem~\ref{thm:Intro.CSZ}.

In what follows, given a Banach space $(\sE, \|\cdot \|_\sE)$, we denote by $C^{(m)}(\R;\sE)$, $m\in \R$, the space of continuous maps $F:\R \rightarrow \sE$ for which there is a constant $C>0$ such that
\begin{equation*}
 \|F(t)\|_{\sE}\leq C \left(1+|t|\right)^m \qquad \forall t\in \R. 
\end{equation*}
This is a Banach space with respect to the norm,
\begin{equation*}
 \left\| F\right\|^{(m)}_\sE:=\sup_{t\in \R} (1+|t|)^{-m} \|F(t)\|_\sE, \qquad F(t) \in C^{(m)}(\R;\sE).  
\end{equation*}

\begin{lemma}[compare~\cite{CSZ:MS17}]\label{lem:CSZ.X1itFX2it}
 Suppose that $X$ and $Y$ are positive compact operators on $\sH$.  
 \begin{enumerate}
 \item[(i)] If $F(t)\in C^{(m)}(\R;\sK)$, $m<-1$, then $X^{it}F(t)Y^{-it}$, $t\in \R$,  is an $L^1$-family in $\sK$.  
 
 \item[(ii)] We define a continuous linear operator $\Theta_{X,Y}:C^{(m)}(\R;\sK)\rightarrow \sK$ by
\begin{equation}\label{eq:CSZ.ThetaXY}
 \Theta_{X,Y}(F) = \int_\R X^{it}F(t)Y^{-it}dt, \qquad F(t)\in C^{(m)}(\R;\sK), 
\end{equation}
where the integral converges in $\sK$ in the Bochner's sense. 
\end{enumerate}
\end{lemma}
\begin{proof}
Let $F(t)\in C^{(m)}(\R;\sK)$, $m<-1$. We know from Lemma~\ref{lem:Com.powers} that the families $X^{it}$, $t\in \R$, and $Y^{-it}$, $t\in \R$, are bounded strongly continuous in $\sL(\sH)$.   Thus, $X^{it}F(t)Y^{-it}$, $t\in \R$, is a bounded strongly continuous family in $\sK$. In particular, it is continuous with respect to the weak operator topology, i.e., it is continuous with respect to the topology defined by the semi-norms $T\rightarrow \scal{T\xi}{\eta}$, $\xi, \eta\in \sH$. Equivalently, for every finite rank operator $R$, the map $t\rightarrow \Tr[RX^{it}F(t)Y^{-it}]$ is continuous on $\R$. Here $\sK'\simeq \sL_1$ and finite-rank operators are dense in $\sL_1$.  
Combining this with the boundedness of the family $X^{it}F(t)Y^{-it}$, $t\in \R$, we deduce this family is continuous with respect to the weak topology of $\sK$, and hence it is weakly measurable. As $\sK$ is separable, Pettis theorem (see, e.g., \cite{HP:AMS74}) then ensures that the family is measurable with respect to the norm topology of $\sK$.   

Moreover, as $m<-1$ and $\|X^{it}\|=\|Y^{-it}\|=1$, we have 
\begin{equation}
 \int_\R \left\|X^{it}F(t)Y^{-it}\right\| dt \leq  \int_\R \left\|F(t)\right\| dt \leq C\|F\|^{(m)}, 
\end{equation}
where we have set $C:= \int (1+|t|)^{m}dt<\infty$. This shows that $X^{it}F(t)Y^{-it}$, $t\in \R$,  is an $L^1$-family in $\sK$ in Bochner's sense. The above inequality further ensures that
\begin{equation*}
 \left\|\int_\R X^{it}F(t)Y^{-it} dt\right\| \leq  \int_\R \left\|X^{it}F(t)Y^{-it}\right\| dt \leq  C\|F\|^{(m)} \qquad \forall F\in C^{(m)}(\R;\sK). 
\end{equation*}
This shows that~(\ref{eq:CSZ.ThetaXY}) defines a bounded operator $\Theta_{X,Y}:C^{(m)}(\R;\sK)\rightarrow \sK$. The proof is complete. 
\end{proof}

\begin{remark}
 A little extra work shows that $\Theta_{X,Y}$ induces a bounded operator from $C^{(m)}(\R;\sL_1)$ to $\sL_1$, and so by interpolation it induces bounded operators from $C^{(m)}(\R;\sL_{q,\infty})$ to $\sL_{q,\infty}$ for all $q>1$. 
\end{remark}

 In what follows for $\Re s>1$ we let $g_s:\R\rightarrow \C$ be the function defined by
\begin{equation}\label{eq:CSZ.gs}
 g_s(t) = -\frac{1}{2} \frac{\sinh \left[ (s-2)t/2\right]}{\sinh(t/2) \cosh\left[(s-1)t/2\right]}, \quad t\neq 0, \qquad g_s(0)= -\frac{1}{2}(s-2). 
\end{equation}
  
\begin{lemma}[compare~\cite{MSZ:LMP22, SZ:Ast23}]\label{lem:CSZ.gs}
 The family $g_s(t)$, $\Re s>1$, is a holomorphic family in $\sS(\R)$. 
\end{lemma}
\begin{proof}
 Set $\Omega=\{\Re s>1\}$. Note that each function $g_s(t)$, $s\in \Omega$, is even. In addition, let $h:\C\rightarrow \C$ be the function defined by
 \begin{equation*}
 h(z)= \frac{\sinh (z)}{z}, \quad z\neq 0, \qquad h(0)=1. 
\end{equation*}
This is an entire function on $\C$ which is~$>0$ on $\R$ and has $i\pi \Z$ as zero-set.  Moreover, we have
\begin{equation*}
 g_s(t)= - \frac{1}{2}(s-2) \frac{h \left[ (s-2)t/2\right]} {h(t/2) \cosh\left[(s-1)t/2\right]}, \qquad \Re s>1, \ t \in \R. 
\end{equation*}
As $\cosh(z)$ is a non-vanishing holomorphic function on the half-plane $\Re z>0$, it follows that $(s,t)\rightarrow g_s(t)$ is a $C^\infty$-function on $\Omega \times \R$ which is holomorphic with respect to $s$. This ensures that $g_s(t)$, $s  \in \Omega$, is at least a holomorphic family in $C^\infty(\R)$. 

In view of~(\ref{eq:CSZ.gs}), for $t>0$ and $s\in \Omega$, we have
\begin{equation*}
 g_s(t)  = \frac{ e^{-(s-2)t/2} -  e^{(s-2)t/2}}{e^{st/2}(1-e^{-t})(1+e^{-(s-1)t})} =  \frac{ e^{-(s-1)t} -  e^{-t}}{1-e^{-t}+e^{-(s-1)t}-e^{-st}}
\end{equation*}
An induction then shows that, for $\ell =0,1, \ldots $ we may write 
\begin{equation*}
 \partial_t^\ell g_s(t)  = \frac{e^{-t}u_s^{(\ell)}(t)+ e^{-(s-1)t}v_s^{(\ell)}(t)}{(1+f_s(t))^{\ell+1}}, \qquad t>0, \ s\in \Omega,
\end{equation*}
where we have set $f_s(t)=e^{-(s-1)t}-e^{-t}-e^{-st}$ and $u_s^{(\ell)}(t)$ and $v_s^{(\ell)}(t)$ are polynomials in $s$, $e^{-t}$ and $e^{-(s-1)t}$. It follows from this (and the fact that $g_s(t)$ is an even function) that, given any $N\geq 1$ and compact $K\subseteq \Omega$, there is a constant $C_{\ell KN}>0$ such that
\begin{equation*}
 \left| \partial_t^\ell g_s(t)\right| \leq C_{\ell KN} |t|^{-N} \qquad \forall (s,t) \in K\times ( \R\setminus (-1,1)). 
\end{equation*}
These estimates together with the fact that $g_s(t)$, $s  \in \Omega$, is a holomorphic family in $C^\infty(\R)$ ensure that $g_s(t)$, $s\in \Omega$, is a holomorphic family in $\sS(\R)$. The proof is complete. 
\end{proof}

\begin{remark}\label{rmk:CSZ.gs}
 As $g_s(t)$, $\Re s>1$, is a Schwartz-class function, we may define its Fourier transform by
 \begin{equation}
\hat{g}_s(t)= \frac{1}{2\pi}\int_{-\infty}^\infty g_s(\xi)e^{-it\xi}d\xi, \qquad t\in \R. 
\end{equation}
The family $\hat{g}_s(t)$, $\Re s>1$, then is a holomorphic family in $\sS(\R)$.  
\end{remark}

\subsection{Proof of the Enhanced CSZ Lemma}\label{subsec:CSZ.proof}
Upon replacing $\sA$ by the algebra $\sA_1$ of Proposition~\ref{prop:App.hol-extension}, which does not change the operator $D$, we may assume that $\sA$ is closed under holomorphic functional calculus. In particular, $a^z\in \sA$ for all $a\in \sA_{++}$ and $z\in \C$. 

Let $a\in {\sA}_{++}$ and $\alpha\in (0,1]$. As in the statement of Theorem~\ref{thm:Intro.CSZhol}, for $\Re s>1$, we set
\begin{equation*}
 T(s):=  \big( a^{\frac12}|D|^{-\alpha} a^{\frac12}\big)^{s} - |D|^{-\alpha s} a^{s}.
\end{equation*}
We shall prove that, for $\Re s>1$, we have
\begin{equation}\label{eq:Com.refined-CSZ1-Ts}
 T(s)\in \sL_{(\alpha \Re s+1-\epsilon)^{-1}p,\infty} \qquad \forall \epsilon \in (0,1).
\end{equation}
Specializing this to $s\in (1,\infty)$ gives Theorem~\ref{thm:Intro.CSZ}. As we shall see, an elaboration of the arguments of the proof will lead to Theorem~\ref{thm:Intro.CSZhol} (see Section~\ref{subsec:CSZ.hol}). Corollary~\ref{cor:Intro.zeta} will then be a by-product of that result (see Section~\ref{subsec:CSZ.zeta}).

Set $X=|D|^{-\alpha}$ and $Y=a^{1/2}|D|^{-\alpha} a^{1/2}$. As mentioned in the Introduction, our approach relies on the integral representation for differences 
$(A^{1/2}BA^{1/2})^z-B^zA^z$, $\Re z>1$,  established by Connes--Sukochev--Zanin~\cite[Lemma~5.2]{CSZ:MS17} and Sukochev--Zanin~\cite[Theorem~5.2.1]{SZ:Ast23}. Applying this formula to $A=a$ and $B=X$ gives
\begin{equation}\label{eq:CSZ.Ts}
 T(s)= \int_\R \hat{g}_s(t)X^{it}F_s(t)Y^{-it}dt-F_s(0), \qquad \Re s>1, 
\end{equation}
where $F_s(t)=F_{1,s}(t)+F_{2,s}(t)$, with 
\begin{equation*}
 F_{1,s}(t):= X^{s-1}\big[X,a^{s-\frac12+it}\big] a^{\frac12},  \qquad F_{2,s}(t):=\big[X,a^{\frac12+it}\big] a^{\frac12} Y^{s-1}.
\end{equation*}
In~\cite{CSZ:MS17, SZ:Ast23} the integral in~(\ref{eq:CSZ.Ts}) is meant in the sense of the weak operator topology. However, the arguments below show that the integral actually converges in $\sK$ in the Bochner's sense. 

Bearing this in mind, let $\epsilon \in (0,1)$. Applying Theorem~\ref{thm:Intro.factorization} for $\beta=1-\epsilon$ and $\gamma=\alpha$ shows there is a bounded operator $\Phi^{(\alpha)}_\epsilon=\Phi^{\alpha}_{1-\epsilon,\alpha}:\sL(\sH)\rightarrow \sK$ such that, for all $b\in \sA$, we have
\begin{align}
[X,b] &= \big[|D|^{-\alpha},b\big] \nonumber \\
& =  |D|^{-\beta}\Phi^{(\alpha)}_\epsilon\left([D,b]\right)|D|^{-\alpha}
                                                  +|D|^{-(1+\alpha)}F[D,b]\Pi_0 +\Pi_0[D,b] F|D|^{-(1+\alpha)},   
\label{eq:CSZ.Dab0}
\end{align}
where $F=D|D|^{-1}$ is the sign of $D$ and  $\Pi_0$ is the orthogonal projection onto $\ker D$. More precisely, in the notation of~(\ref{eq:Com.Phimbgam1}) we have 
\begin{equation*}
 \Phi^{(\alpha)}_\epsilon(T) =-c(\alpha)F\Psi^{-\alpha}_{2-\epsilon,\alpha}(T)-c(\alpha)\Psi^{-\alpha}_{1-\epsilon,\alpha+1}(T)F, \qquad T\in \sL(\sH),
\end{equation*}
where $\Psi^{-\alpha}_{2-\epsilon,\alpha}(T)$ and $\Psi^{-\alpha}_{1-\epsilon,\alpha+1}(T)$ are DOIs given by~(\ref{eq:Com.Psiadd'}). As $|D|^{-1}=X^{1/\alpha}$ and $ |D|^{-\alpha}=a^{-1/2}Ya^{-\frac12}$, we may rewrite the above formula in the form,
\begin{align}
[X,b] =  X^{\beta/\alpha}\Phi^{(\alpha)}_\epsilon\left([D,b]\right)a^{-1/2}Ya^{-\frac12}
                                                  +X^{1+(1/\alpha)}F[D,b]\Pi_0 + \Pi_0[D,b] FX^{1+(1/\alpha)}.   
\label{eq:CSZ.Dab}
\end{align}

By using~(\ref{eq:CSZ.Dab}) and the vanishing of $X^{s-1}$ on $\ran \Pi_0=\ker D$ we get
\begin{equation}\label{eq:CSZ.Fst1}
 F_{1,s}(t)=  X^{(\beta/\alpha)+s-1}\Phi^{(\alpha)}_\epsilon\big([D,a^{s-\frac12 +it}]\big)a^{-\frac12} Y + X^{(1/\alpha)+s}F[D,a^{s-\frac12 +it}]\Pi_0a^{\frac12}. 
\end{equation}
It is convenient to rewrite this in the form, 
\begin{equation}\label{eq:CSZ.F1s}
  F_{1,s}(t)=  X^{(\beta/\alpha)+s-1}G_{1,s}(t) Y + X^{(1/\alpha)+s}H_{1,s}(t), 
\end{equation}
where we have set
\begin{equation*}
 G_{1,s}(t)=\Phi^{(\alpha)}_\epsilon\big([D,a^{s-\frac12 +it}]\big)a^{-\frac12}, \qquad H_{1,s}(t)= F[D,a^{s-\frac12 +it}]\Pi_0a^{\frac12}. 
\end{equation*}

Similarly, we obtain 
\begin{align}
 F_{2,s}(t) = & \  X^{(\beta/\alpha)}\Phi^{(\alpha)}_\epsilon\big([D,a^{\frac12 +it}]\big)a^{-\frac12} Y^s + X^{1+(1/\alpha)}F[D,a^{\frac12 +it}]\Pi_0a^{\frac12}Y^{s-1} \nonumber\\
 & + \Pi_0[D,a^{\frac12 +it}] FX^{1+(1/\alpha)}a^{\frac12}Y^{s-1}. 
 \label{eq:CSZ.Fst2}
 \end{align}
Equivalently, 
\begin{equation}\label{eq:CSZ.F2s}
 F_{2,s}(t)=\tilde{F}_{2,s}(t)+\Pi_0 R_s(t)Y^{s-1}, \qquad \tilde{F}_{2,s}(t):= X^{(\beta/\alpha)} G_2(t)Y^s+  X^{1+(1/\alpha)}H_2(t)Y^{s-1},  
\end{equation}
where we have set 
\begin{equation*}
 G_2(t)= \Phi^{(\alpha)}_\epsilon\big([D,a^{\frac12 +it}]\big)a^{-\frac12}, \quad H_2(t)=F[D,a^{\frac12 +it}]\Pi_0a^{\frac12}, \quad 
 R_s(t)=[D,a^{\frac12 +it}] FX^{1+(1/\alpha)}a^{\frac12}. 
\end{equation*}
Note that as $X^{it}\Pi_0=0$, we have
\begin{equation}\label{eq:CSZ.XF2Y}
 X^{it}F_{2,s}(t)Y^{-it}= X^{it}\tilde{F}_{2,s}(t)Y^{-it}. 
\end{equation}

Given any $z\in \C$, Lemma~\ref{lem:CSZ.Daz} implies that $[D,a^{z+it}]$, $t\in \R$, is a family in $C^{(1)}(\R;\sL(\sH))$. As $\Pi_0$ is bounded and has finite-rank, and hence is compact, we deduce that $H_{1,s}(t)$, $t\in \R$,  and $H_2(t)$, $t\in \R$,  are families in  $C^{(1)}(\R;\sK)$. As $\hat{g}_s(t)$ is in $\sS(\R)$, we further see that 
$\hat{g}_s(t)H_{1,s}(t)$ and $\hat{g}_s(t)H_{2}(t)$ are in $C^{(m)}(\R;\sK)$ for any $m<-1$. Lemma~\ref{lem:CSZ.X1itFX2it} then ensures that 
$\hat{g}_s(t)X^{it}H_{1,s}(t)Y^{-it}$ and $\hat{g}_s(t)X^{it}H_{2}(t)Y^{-it}$ are in $L^1(\R;\sK)$. It follows that  
$\hat{g}_s(t)X^{it}X^{(1/\alpha)+s}H_{1,s}(t)Y^{-it}$ is in $L^1(\R;\sK)$ and, using the notation of Lemma~\ref{lem:CSZ.X1itFX2it}, we have
\begin{equation}\label{eq:CSZ.Theta-H1s}
 \int_\R\hat{g}_s(t)X^{it}\big(X^{(1/\alpha)+s}H_{1,s}(t)\big)Y^{-it}dt= X^{(1/\alpha)+s} \Theta_{X,Y}\big[\hat{g}_sH_{1,s}\big]. 
\end{equation}
Likewise,  $\hat{g}_s(t)X^{it}X^{1+(1/\alpha)}H_{2}(t)Y^{s-1}Y^{-it}$ is in $L^1(\R;\sK)$, and we have
\begin{equation}\label{eq:CSZ.Theta-H2}
 \int_\R\hat{g}_s(t)X^{it}\big(X^{1+(1/\alpha)}H_{2}(t)Y^{s-1}\big)Y^{-it}dt= X^{1+(1/\alpha)} \Theta_{X,Y}\big[\hat{g}_sH_{2}\big]Y^{s-1}. 
\end{equation}

In addition, as $\Phi^{(\alpha)}_\epsilon$ is a bounded operator from $\sL(\sH)$ to $\sK$, we see that $\Phi^{(\alpha)}_\epsilon ([D,a^{z+it}])$ is a family in $C^{(1)}(\R;\sK)$, and hence $G_{1,s}(t)$ and $G_2(t)$ are families in $C^{(1)}(\R;\sK)$. As above this implies that 
$\hat{g}_s(t)X^{it}X^{(\beta/\alpha)+s-1}G_{1,s}(t)YY^{-it}$ is in $L^1(\R;\sK)$, and we have
\begin{equation}\label{eq:CSZ.Theta-G1s}
 \int_\R\hat{g}_s(t)X^{it}\big(X^{(\beta/\alpha)+s-1}G_{1,s}(t)Y\big)Y^{-it}dt= X^{(\beta/\alpha)+s-1} \Theta_{X,Y}\big[\hat{g}_sG_{1,s}\big]Y. 
\end{equation}
Similarly, $\hat{g}_s(t)X^{it}X^{\beta/\alpha}G_{2}(t)Y^{s}Y^{-it}$ is in $L^1(\R;\sK)$, and we have
\begin{equation}\label{eq:CSZ.Theta-G2}
 \int_\R\hat{g}_s(t)X^{it}\big(X^{\beta/\alpha}G_{2}(t)Y^{s}\big)Y^{-it}dt=X^{\beta/\alpha} \Theta_{X,Y}\big[\hat{g}_sG_{2}\big]Y^s. 
\end{equation}

Combining~(\ref{eq:CSZ.F1s}) with~(\ref{eq:CSZ.Theta-H1s}) and~(\ref{eq:CSZ.Theta-G1s}) gives
\begin{equation*}
 \int_\R \hat{g}_s(t)X^{it}F_{1,s}(t)Y^{-it}dt = X^{(\beta/\alpha)+s-1} \Theta_{X,Y}\big[\hat{g}_s G_{1,s}\big]Y + X^{(1/\alpha)+s} \Theta_{X,Y}\big[\hat{g}_s H_{1,s}\big]. 
\end{equation*}
Similarly, by combining~(\ref{eq:CSZ.F2s})--(\ref{eq:CSZ.XF2Y}) with~(\ref{eq:CSZ.Theta-H2}) and~(\ref{eq:CSZ.Theta-G2}) we get
\begin{align*}
 \int_\R \hat{g}_s(t)X^{it}F_{2,s}(t)Y^{-it}dt& =  \int_\R \hat{g}_s(t)X^{it}\tilde{F}_s(t)Y^{-it}dt\\
 & =  X^{\beta/\alpha} \Theta_{X,Y}\big[\hat{g}_s G_{2}\big]Y^s + X^{1+(1/\alpha)} \Theta_{X,Y}\big[\hat{g}_s H_{2}\big]Y^{s-1}. 
\end{align*}
As $F_s(t)=F_{1,s}(t)+F_{2,s}(t)$, we then obtain
\begin{equation}\label{eq:CSZ.ThetaF}
  \int_\R \hat{g}_s(t)X^{it}F_{s}(t)Y^{-it}dt = X^{\beta/\alpha}A(s)Y +  X^{1+(1/\alpha)}B(s),
\end{equation}
where we have set
\begin{equation*}
 A(s)= X^{s-1} \Theta_{X,Y}\big[\hat{g}_s G_{1,s}\big]+\Theta_{X,Y}\big[\hat{g}_s G_{2}\big]Y^{s-1}, \quad 
 B(s)=X^{s-1} \Theta_{X,Y}\big[\hat{g}_s H_{1,s}\big]+\Theta_{X,Y}\big[\hat{g}_s H_{2}\big]Y^{s-1}. 
\end{equation*}

Moreover, setting $t=0$ in (\ref{eq:CSZ.F1s}) and~(\ref{eq:CSZ.F2s}) gives
\begin{gather*}
  F_{1,s}(0)=  X^{(\beta/\alpha)+s-1}G_{1,s}(0) Y + X^{(1/\alpha)+s}H_{1,s}(0),\\
  F_{2,s}(0)= X^{(\beta/\alpha)} G_2(0)Y^s+  X^{1+(1/\alpha)}H_2(0)Y^{s-1} + \Pi_0 R_s(0)Y^{s-1}. 
\end{gather*}
Thus, 
\begin{align}
 F_s(0)=F_{1,s}(0)+F_{2,s}(0) &=X^{\beta/\alpha}\left( X^{s-1}G_{1,s}(0)+G_2(0)Y^{s-1}\right) Y \nonumber \\ 
 & + 
 X^{1+(1/\alpha)}\left( X^{s-1}H_{1,s}(0)+H_2(0)Y^{s-1}\right) + \Pi_0 R_s(0)Y^{s-1}. 
 \label{eq:CSZ.F(0)}
\end{align}

For $\Re s>1$, let $\tilde{\Theta}_{X,Y}^{(s)}:C^{(1)}(\R;\sK) \rightarrow \sK$ be the operator defined by
\begin{equation}\label{eq:CSZ.tThetaXY}
 \tilde{\Theta}_{X,Y}^{(s)}[K]= {\Theta}_{X,Y}[\hat{g}_s K]-K(0), \qquad K(t)\in C^{(1)}(\R;\sK). 
\end{equation}
Combining~(\ref{eq:CSZ.Ts}) with~(\ref{eq:CSZ.ThetaF}) and~(\ref{eq:CSZ.F(0)}) gives 
\begin{align}
 T(s) & = \int_\R \hat{g}_s(t)X^{it}F_{s}(t)Y^{-it}dt -F_s(0)\nonumber\\
 &  = X^{\beta/\alpha}\tilde{A}(s)Y +  X^{1+(1/\alpha)}\tilde{B}(s) - \Pi_0 R(s), 
 \label{eq:CSZ.Ts-As-Bs-Rs}
 \end{align}
where we have set
\begin{gather}\label{eq:CSZ.tAs}
 \tilde{A}(s):=X^{s-1} \tilde{\Theta}_{X,Y}^{(s)}\big[ G_{1,s}\big]+\tilde{\Theta}_{X,Y}^{(s)}\big[ G_{2}\big]Y^{s-1},  \\
 \tilde{B}(s):=X^{s-1} \tilde{\Theta}_{X,Y}^{(s)}\big[H_{1,s}\big]+\tilde{\Theta}_{X,Y}^{(s)}\big[ H_{2}\big]Y^{s-1},
 \label{eq:CSZ.tBs}\\
 R(s):=R_s(0)Y^{s-1}= [D,a^{\frac12}] FX^{1+(1/\alpha)}a^{\frac12}Y^{s-1}.
\end{gather}

Set $c=\Re s$. The operators  $\tilde{\Theta}_{X,Y}^{(s)}[G_{1,s}]$ and $\tilde{\Theta}_{X,Y}^{(s)}\big[ G_{2}\big]$ are both bounded. Likewise, 
$\tilde{\Theta}_{X,Y}^{(s)} [H_{1,s}]$ and $\tilde{\Theta}_{X,Y}^{(s)}\big[ H_{2}\big]$ are bounded operators. Moreover, as $X$ and $Y$ are operators in $\sL_{\alpha^{-1}p,\infty}$, the powers $X^{s-1}$ and $Y^{s-1}$ are both operators in $\sL_{(c\alpha-\alpha)^{-1}p,\infty}$. It then follows that $\tilde{A}(s)$ and $\tilde{B}(s)$ are operators in $\sL_{(c\alpha-\alpha)^{-1}p,\infty}$. 

The fact that $X$ is in $\sL_{\alpha^{-1}p,\infty}$ further implies that $X^{\beta/\alpha}\in \sL_{\beta^{-1}p,\infty}$. As
 $\tilde{A}(s)\in  \sL_{(c\alpha-\alpha)^{-1}p,\infty}$ and $Y\in \sL_{\alpha^{-1}p,\infty}$, by using H\"older's inequality (Proposition~\ref{prop:Hoelder}) we then get
\begin{equation}\label{eq:CSZ.As-sL}
 X^{\beta/\alpha}\tilde{A}(s)Y \in \sL_{(c\alpha+\beta)^{-1}p,\infty}. 
\end{equation}
Similarly, the operator $X^{1+(1/\alpha)}$ is in $\sL_{(\alpha+1)^{-1}p,\infty}$. As $\tilde{B}(s)\in \sL_{(c\alpha-\alpha)^{-1}p,\infty}$, by using H\"older's inequality we get 
\begin{equation}\label{eq:CSZ.Bs-sL}
 X^{1+(1/\alpha)}\tilde{B}(s) \in \sL_{(c\alpha+1)^{-1}p,\infty}\subseteq \sL_{(c\alpha+\beta)^{-1}p,\infty}. 
\end{equation}
In addition, the operator $\Pi_0 R(s)$ is a finite-rank bounded operator, and hence is contained in $\sL_{(c\alpha+\beta)^{-1}p,\infty}$. Combining this with~(\ref{eq:CSZ.Ts-As-Bs-Rs}) and 
(\ref{eq:CSZ.As-sL})--(\ref{eq:CSZ.Bs-sL}) we then see that
\begin{equation*}
 T(s)\in  \sL_{(c\alpha+\beta)^{-1}p,\infty}. 
\end{equation*}
 As $\beta=1-\epsilon$ this gives~(\ref{eq:Com.refined-CSZ1-Ts}). This completes the proof of Theorem~\ref{thm:Intro.CSZ}.  \hfill $\Box$ 
 
 \begin{remark}\label{rmk:CSZ.alpha}
Eq.~(\ref{eq:CSZ.Fst1}) and Eq.~(\ref{eq:CSZ.Fst2}) are two of the main steps in the proof of Theorem~\ref{thm:Intro.CSZ} above. For getting these results it was crucial to be able to factor out $|D|^{-\alpha}$ in the factorization~(\ref{eq:CSZ.Dab0}) for $[|D|^{-\alpha},b]$. 
\end{remark}

\subsection{Holomorphic version}\label{subsec:CSZ.hol}
An elaboration of the arguments of the proof of Theorem~\ref{thm:Intro.CSZ} leads to its holomorphic version, Theorem~\ref{thm:Intro.CSZhol}. In what follows we keep on using the notation of Section~\ref{subsec:CSZ.proof}.

Prior to proving Theorem~\ref{thm:Intro.CSZhol}, we will need the following lemma.

\begin{lemma} The following holds. 
\begin{enumerate}
 \item[(i)]
 $\tilde{\Theta}_{X,Y}^{(s)}[ G_{1,s}]$, $\Re s>1$, and $\tilde{\Theta}_{X,Y}^{(s)}[G_2]$, $\Re s>1$, are
  holomorphic families of compact operators. 
 
  \item[(ii)] The same holds for the families $\tilde{\Theta}_{X,Y}^{(s)}[ H_{1,s}]$, $\Re s>1$, and $\tilde{\Theta}_{X,Y}^{(s)}[ H_{2}]$, $\Re s>1$,
\end{enumerate}
\end{lemma}
\begin{proof}
Suppose that $K_s$, $\Re s>1$, is a holomorphic family in $C^{(1)}(\R; \sK)$. We know from Lemma~\ref{lem:CSZ.gs} and Remark~\ref{rmk:CSZ.gs} that 
$\hat{g}_s(t)$, $\Re s>1$, is a holomorphic family in $\sS(\R)$. 
Therefore, $\hat{g}_s K_s$, $\Re s>1$, is a holomorphic family in $C^{(m)}(\R;\sK)$ for any $m<-1$.  
As Lemma~\ref{lem:CSZ.X1itFX2it} asserts that $\Theta_{X,Y}:C^{(m)}(\R;\sK)\rightarrow \sK$ is a continuous linear operator, we see that $\Theta_{X,Y}[\hat{g}_s K_s]$, $\Re s>1$, is a holomorphic family in $\sK$. 
Clearly, $K_s(0)$, $\Re s>1$, is a holomorphic family in $\sK$, and so in view of~(\ref{eq:CSZ.tThetaXY}) we see that $\tilde{\Theta}_{X,Y}^{(s)}[ K_{s}]$, $\Re s>1$, is a holomorphic family of compact operators. 

Applying the above observation to $K_s=G_2$ and $K_s=H_2$ shows that $\tilde{\Theta}_{X,Y}^{(s)}[ G_{2}]$, $\Re s>1$, and $\tilde{\Theta}_{X,Y}^{(s)}[ H_{2}]$, $\Re s>1$, are holomorphic families of compact operators. Moreover, in order to prove the claim for $\tilde{\Theta}_{X,Y}^{(s)}[ G_{1,s}]$ and $\tilde{\Theta}_{X,Y}^{(s)}[ H_{1,s}]$ it is enough to show that $G_{1,s}$, $\Re s>1$, and $H_{1,s}$, $\Re s>1$, are holomorphic families in $C^{(1)}(\R; \sK)$. 

Set $\Omega=\{\Re s>1\}$.  Lemma~\ref{lem:CSZ.Daz} asserts that $[D,a^z]$, $z\in \C$, is a holomorphic family in $\sL(\sH)$, and, for any (closed) vertical strip $\Sigma\subseteq \C$, there is $C_{\Sigma}(a)>0$ such that
\begin{equation*}
 \big\|[D,a^z]\big\| \leq C_{\Sigma}(a) \left( 1+|\Im z|\right) \qquad \forall z\in \Sigma.  
\end{equation*}
 Therefore, if we set $V_s(t)=[D,a^{s-1/2+it}]$, $\Re s>1$, $t\in \R$, then $V_s(t)$, $\Re s>1$, is a holomorphic family in $C(\R;\sL(\sH))$, and, for any compact $L\subseteq \Omega$, there is $C_{L}(a)>0$ such that 
\begin{equation*}
  \big\|V_{s}(t)\big\| \leq C_L(a) (1+|t|)\qquad \forall (s,t)\in L\times \R. 
\end{equation*}
This ensures that $V_s(t)$, $\Re s>1$, is a holomorphic family in $C^{(1)}(\R;\sL(\sH))$. 

By definition $H_{1,s}(t)=FV_{s}(t)\Pi_0a^{\frac12}$. Here $\Pi_0$ is a bounded finite-rank operator, and hence is compact, and so we see that $H_{1,s}$, $\Re s>1$,  is a holomorphic family in $C^{(1)}(\R;\sK)$. Moreover, as 
$\Phi_\epsilon^{(\alpha)}$ is a continuous linear operator from $\sL(\sH)$ to $\sK$ and $G_{1,s}=\Phi^{(\alpha)}_\epsilon[V_s(t)]a^{-\frac12}$, we see that 
$G_{1,s}$, $\Re s>1$, is a holomorphic family in $C^{(1)}(\R;\sK)$ as well. This completes the proof.   
\end{proof}
  
 We are now in a position to prove Theorem~\ref{thm:Intro.CSZhol}.
  \begin{proof}[Proof of Theorem~\ref{thm:Intro.CSZhol}] 
  Let $c> 1$. As $X$ is a positive operator in $\sL_{\alpha^{-1}p,\infty}$, Lemma~\ref{lem:Tauberian.powers} ensures that $X^z$, $\Re z>c-1$, is a holomorphic family in $\sL_{(c\alpha-\alpha)^{-1}p,\infty}$, i.e., $X^{s-1}$, $\Re s>c$, is a holomorphic family in $\sL_{(c\alpha-\alpha)^{-1}p,\infty}$. 
 Likewise, $Y^{s-1}$, $\Re s>c$, is holomorphic family in $\sL_{(c\alpha-\alpha)^{-1}p,\infty}$. Combining this with~(\ref{eq:CSZ.tAs})--(\ref{eq:CSZ.tBs}) and the claim above shows that 
 \begin{equation}\label{eq:CSZ.tAs-tBs-Hol}
 \tilde{A}(s) \in \Hol\left( \Re s>c; \sL_{(c\alpha-\alpha)^{-1}p,\infty}\right) \quad \textup{and} \quad 
 \tilde{B}(s) \in \Hol\left( \Re s>c; \sL_{(c\alpha-\alpha)^{-1}p,\infty}\right).  
\end{equation}
 
 If $r_1^{-1}+r_2^{-1}=r^{-1}$, then H\"older's inequality implies we have a continuous bilinear map from $\sL_{r_1,\infty}\times \sL_{r_2,\infty}$ to $\sL_{r,\infty}$. Thus, if $Z_1\in \sL_{r_1,\infty}$ and $Z_2(s)$, $\Re s>c$, is a holomorphic family in $\sL_{r_2,\infty}$, then $Z_1Z_2(s)$ is a holomorphic family in $\sL_{r,\infty}$. Using this remark and the fact that $X^{\beta/\alpha}\in \sL_{\beta^{-1}p,\infty}$ and $Y\in \sL_{\alpha^{-1}p,\infty}$ together with~(\ref{eq:CSZ.tAs-tBs-Hol}) shows that
 \begin{equation}\label{eq:CSZ.XAsY-Hol}
X^{\beta/\alpha} \tilde{A}(s)Y \in \Hol\left( \Re s>c; \sL_{(c\alpha+\beta)^{-1}p,\infty}\right). 
\end{equation}
 Similarly, as $X^{1+(1/\alpha)}\in \sL_{(\alpha+1)^{-1}p,\infty}$ we see that
 \begin{equation}
 X^{1+(1/\alpha)}\tilde{B}(s) \in \Hol\left( \Re s>c; \sL_{(c\alpha+1)^{-1}p,\infty}\right). 
\end{equation}
 Likewise, $R(s)=[D,a^{\frac12}] FX^{1+(1/\alpha)}a^{\frac12}Y^{s-1}$, $\Re s>c$, is a holomorphic family in $\sL_{(c\alpha+1)^{-1}p,\infty}$, and hence 
 \begin{equation}\label{eq:CSZ.PiRs-Hol}
 \Pi_0R(s)\in \Hol\left( \Re s>c; \sL_{(c\alpha+1)^{-1}p,\infty}\right).  
\end{equation}
 Combining~(\ref{eq:CSZ.Ts-As-Bs-Rs}) with (\ref{eq:CSZ.XAsY-Hol})--(\ref{eq:CSZ.PiRs-Hol}), and recalling that $\beta=1-\epsilon$, we then obtain
\begin{equation*}
 T(s) \in \Hol\left( \Re s>c; \sL_{(c\alpha+1-\epsilon)^{-1}p,\infty}\right) \qquad \forall \epsilon\in (0,1). 
\end{equation*}

Let $q>0$, and set $c=\max(1,\alpha^{-1}(q^{-1}p-1))$. The inequality $\alpha^{-1}(q^{-1}p-1)\leq c$ ensures that $(c\alpha+1)^{-1}p\leq q$. Let $\epsilon \in (0,1)$. Substituting $c+2\epsilon\alpha^{-1}$ for $c$ shows that  
$T(s)$, $\Re s>c+2\epsilon\alpha^{-1}$, is a holomorphic family in $\sL_{(c\alpha+1+\epsilon)^{-1}p,\infty}$. As $(c\alpha+1+\epsilon)^{-1}p< (c\alpha+1)^{-1}p\leq q$, we have a continuous inclusion of $\sL_{(c\alpha+1+\epsilon)^{-1}p,\infty}$ into 
$\sL_{q}$. This allows us to regard $T(s)$, $\Re s>c+2\epsilon\alpha^{-1}$, as a holomorphic family in 
$ \sL_{q}$. As this is true for all $\epsilon\in (0,1)$, we deduce that  
\begin{equation*}
 T(s) \in \Hol\left( \Re s>c; \sL_q\right). 
\end{equation*}
This proves Theorem~\ref{thm:Intro.CSZhol}.
\end{proof}

 \subsection{Differences of zeta functions}\label{subsec:CSZ.zeta}
In this subsection we prove Corollary~\ref{cor:Intro.zeta}. Let $a\in \sA_{++}$ and $\alpha\in (0,1]$. It follows from Lemma~\ref{lem:Tauberian.powers} that $(a^{\frac{1}{2}}|D|^{-\alpha}a^{\frac{1}{2}})^{s}$ and $|D|^{-\alpha s}a^{s}$ are holomorphic families of trace-class operators for $\Re s>\alpha^{-1}p$. Therefore, the zeta functions $\Tr[(a^{1/2}|D|^{-\alpha}a^{1/2})^{s}]$ and $\Tr[a^{s}|D|^{-\alpha s}]$ are holomorphic functions on the half-plane $\Re s> \alpha^{-1}p$. 

In the special case $q=1$ Theorem~\ref{thm:Intro.CSZhol} shows that, if we set $c=\max(1,\alpha^{-1}(p-1))$, then we have
\begin{equation}
  \left( a^{\frac12}|D|^{-\alpha} a^{\frac12}\right)^{s} - |D|^{-\alpha s} a^{s} \in \Hol\left( \Re s>c; \sL_1\right). 
\end{equation}

Assume now that $\delta:=\min(1,p-\alpha)>0$. We then have
\begin{equation*}
 c=\alpha^{-1}p + \max\left(1-\alpha^{-1}p,-\alpha^{-1}\right)= \alpha^{-1}p - \alpha^{-1}\min\left(p-\alpha,1\right)= \alpha^{-1}(p-\delta). 
\end{equation*}
Thus, 
\begin{equation}
  \left( a^{\frac12}|D|^{-\alpha} a^{\frac12}\right)^{s} - |D|^{-\alpha s} a^{s} \in \Hol\left( \Re s>\alpha^{-1}(p-\delta); \sL_1\right).  
\end{equation}
 It then follows that the function $s\rightarrow \Tr[(a^{1/2}|D|^{-\alpha}a^{1/2})^{s}]-\Tr[a^{s}|D|^{-\alpha s}]$ is holomorphic on the half-plane $\Re s>\alpha^{-1}(p-\delta)$. This completes the proof of Corollary~\ref{cor:Intro.zeta}.\hfill $\Box$ 

\section{Tauberian Theorem}\label{chap:Tauberian}

This section is devoted to the proof of the Tauberian Theorem (Theorem~\ref{thm:Intro.Tauberian}). As a byproduct of the approach of this section, we will even get a stronger version for spectral triples of the Tauberian theorem of~\cite{MSZ:LMP22} without any extra regularity or summability conditions (see Remark~\ref{rmk:Tauberian.comparison-MSZ}).

Throughout this section we let $(\sA, \sH,D)$ be a $p$-summable spectral triple, with $p>0$ arbitrary, and we let $\tau:\overline{\sA}\rightarrow \C$ be a non-zero positive linear map.

\subsection{From Condition ($\textup{Z}_0$) to Condition (W)} In our setting, the Tauberian condition of~\cite{MSZ:LMP22} is formulated as follows.

\begin{condition}[Tauberian Condition of \cite{MSZ:LMP22}] \label{cnd:TauberianMSZ} 
 For every $a\in \overline{\sA}_{+}$, the function
 \begin{equation*}
 \Tr\big[a^{s}|D|^{-s}\big] - p\tau\left[a^{p}\right](s-p)^{-1}, \qquad \Re s>p,
\end{equation*}
has a unique continuous extension to the closed half-plane $\Re s \geq p$. 
\end{condition}

It is too strong a requirement to have the condition hold for all positive elements of the $C^*$-algebra $\overline{\sA}$. In fact, in examples where $\sA$ is an algebra of \psidos\ it may be technically cumbersome to check Condition~\ref{cnd:TauberianMSZ} for all elements of $\overline{\sA}_+$, whereas it is comparatively straightforward to check Weyl laws for the operators $aD^2a$, $a\in \sA_{++}$ (i.e., Condition~(W)). 

In addition, as mentioned in the Introduction, it was asked by Alain Connes to relate Condition~\ref{cnd:TauberianMSZ} to more standard Tauberian conditions involving the zeta functions $\Tr[a|D|^{-s}]$, $a\in \sA$, instead of the functions $\Tr[a^s|D|^{-s}]$, $a\in \overline{\sA}_+$. He also pointed out that closedness under holomorphic functional calculus could come into play. 

In fact, what is relevant in our setting is a weakening of Condition~\ref{cnd:TauberianMSZ} provided by Condition~($\textup{Z}_0$) below. For one thing, that condition implies Condition~(W) (see Proposition~\ref{prop:Tauberian.Weyl-W} below). For another thing, it will be shown that Condition~($\textup{Z}_0$) is implied by both Condition~(Z) and Condition~(H) (see Lemma~\ref{lem:Tauberian.Z-Z0} and Lemma~\ref{lem:Tauberian.H-Z} below). This will prove Theorem~\ref{thm:Intro.Tauberian}. 

\begin{condition}[Condition ($\textup{Z}_0$)] 
 There is a holomorphically closed Fr\'echet subalgebra $\sB\subseteq \overline{\sA}$ containing $\sA$ such that, for every $a\in {\sB}_{++}$, the function
 \begin{equation*}
 \Tr\big[a^{s}|D|^{-s}\big] - p\tau\left[a^{p}\right](s-p)^{-1}, \qquad \Re s>p,
\end{equation*}
has a unique continuous extension to the closed half-plane $\Re s \geq p$. 
\end{condition}

\begin{lemma}\label{lem:Tauberian.Z0-W}
Assume Condition~$\textup{(Z}_0\textup{)}$ holds. Then, for every $\alpha\in (0,1]$ such that $\alpha<p$, we have the Weyl laws
\begin{equation}\label{eq:Tauberian.Weyl}
 \lim_{j\rightarrow \infty} j^{\frac{\alpha}{p}}\lambda_j\big(a^{\frac12}|D|^{-\alpha}a^{\frac12}\big)= \tau\big[a^{\frac{p}{\alpha}}\big]^{\frac{\alpha}{p}} \qquad \forall a\in \sA_{++}.
\end{equation}
\end{lemma}
\begin{proof}
 Assume Condition~$\textup{(Z}_0\textup{)}$ holds. Let $\alpha \in (0,1]$ be such that $\alpha<p$, so that $\delta:=\min(1, p-\alpha)>0$, and let $a\in {\sA}_{++}$. Define
  \begin{equation*}
 N(\lambda):= \#\left\{ j; \ \lambda_j\big(a^{\frac12}|D|^{-\alpha} a^{\frac12}\big)>\lambda^{-1}\right\}, \qquad \lambda>0. 
\end{equation*}
Alternatively, $N(\lambda^{-1})$ is the counting function of $a^{1/2}|D|^{-\alpha} a^{1/2}$. In particular, $N(\lambda)$ is a non-negative non-decreasing step function such that $N(\lambda)=0$ for $\lambda<\lambda_0(a^{1/2}|D|^{-\alpha} a^{1/2})^{-1}$. As $dN(\lambda)=\sum \delta_{\lambda_j(a^{1/2}|D|^{-\alpha} a^{1/2})^{-1}}$, for $\Re s>\alpha^{-1}p$ we have
\begin{equation*}
 \Tr\left[ \big(a^{\frac12}|D|^{-\alpha} a^{\frac12}\big)^{s}\right] = 
 \sum_{j\geq 0} \left(\lambda_j\big(a^{\frac12}|D|^{-\alpha} a^{\frac12}\big)^{-1}\right)^{-s}=  \int_0^\infty \lambda^{-s} dN(\lambda). 
\end{equation*}
Moreover, the fact that  $N(\lambda^{-1})$ is the counting function of $a^{1/2}|D|^{-\alpha} a^{1/2}$ further implies that
\begin{equation}\label{eq:Tauberian.counting-eigenvalues}
 \lim_{j\rightarrow \infty} j\lambda_j\big(a^{\frac12}|D|^{-\alpha} a^{\frac12}\big)^{\frac{p}{\alpha}}=
  \lim_{\lambda \rightarrow 0^+} \lambda^{\frac{p}{\alpha}}N(\lambda^{-1})=  
  \lim_{\lambda \rightarrow \infty} \lambda^{-\frac{p}{\alpha}}N(\lambda), 
\end{equation}
provided any of these limits exists.

Bearing this in mind, for $\Re s>\alpha^{-1}p$ set
\begin{gather*}
 F(s)= \Tr\left[ \big(a^{\frac12}|D|^{-\alpha} a^{\frac12}\big)^{s}\right]  -\Tr\big[a^{s}|D|^{-\alpha s}\big], \\
 G(s) = \Tr\big[a^{s}|D|^{-\alpha s}\big]- \alpha^{-1}p\tau\big[a^{\frac{p}{\alpha}}\big](s-\alpha^{-1}p)^{-1}.  
\end{gather*}
We then have
\begin{equation}\label{eq:Tauberian.Tr-FG}
\Tr\left[ \big(a^{\frac12}|D|^{-\alpha} a^{\frac12}\big)^{s}\right] =\alpha^{-1}p\tau\big[a^{\frac{p}{\alpha}}\big](s-\alpha^{-1}p)^{-1} +  F(s) + G(s). 
\end{equation}
 
 By Corollary~\ref{cor:Intro.zeta} the function $F(s)$ has an analytic extension to the half-plane $\Re s>\alpha^{-1}(p-\delta)$ with $\delta$ as above. Moreover, let $\sB$ be as in Condition~$\textup{(Z}_0\textup{)}$. As $\sB$ is holomorphically closed and $a\in \sA_{++}\subseteq \sB_{++}$, we have $a^{1/\alpha}\in \sB_{++}$. Therefore, by assumption 
 $G(s)=\Tr[(a^{1/\alpha})^{ \alpha s }|D|^{-\alpha s}]- p\tau[(a^{1/\alpha})^p](\alpha s-p)^{-1}$ has a continuous extension to the half-plane $\Re s\geq \alpha^{-1}p$. It follows that the function 
 \begin{equation*}
\Tr\left[ \big(a^{\frac12}|D|^{-\alpha} a^{\frac12}\big)^{s}\right] -\alpha^{-1}p\tau\big[a^{\frac{p}{\alpha}}\big](s-\alpha^{-1}p)^{-1}, \qquad \Re s>\alpha^{-1}p,
\end{equation*}
 has a continuous extension to the half-plane $\Re s\geq \alpha^{-1}p$. Combining this with~(\ref{eq:Tauberian.Tr-FG}) and Ikehara's Tauberian theorem (see, e.g.,~\cite[Theorem~14.1]{Sh:Springer01}), and using~(\ref{eq:Tauberian.counting-eigenvalues}), we then obtain
 \begin{equation}\label{eq:Tauberian.Ikehara}
 \lim_{j\rightarrow \infty} j\lambda_j\big(a^{\frac12}|D|^{-\alpha} a^{\frac12}\big)^{\frac{p}{\alpha}}= \lim_{\lambda\rightarrow \infty} \lambda^{-\frac{p}{\alpha}} N(\lambda)= \tau\big[a^{\frac{p}{\alpha}}\big]. 
\end{equation}

Equivalently, we have
\begin{equation*}
 \lim_{j\rightarrow \infty} j^{\frac{\alpha}{p}}\lambda_j\big(a^{\frac12}|D|^{-\alpha} a^{\frac12}\big)= \tau\big[a^{\frac{p}{\alpha}}\big]^{\frac{\alpha}{p}}.
\end{equation*}
That is, the Weyl law~(\ref{eq:Tauberian.Weyl}) holds for $a$. The proof is complete.  
\end{proof}

\begin{remark}\label{rmk:Tauberian.comparison-MSZ}
 The equality~(\ref{eq:Tauberian.Ikehara}) is merely the contents of the Tauberian theorem of~\cite{MSZ:LMP22} in its version for spectral triples (i.e., \cite[Theorem~1.4]{MSZ:LMP22}). In~\cite{MSZ:LMP22} the result is obtained by assuming Lipschitz regularity and $p$-summability with $p>2$. Both assumptions are unnecessary in our approach.   
\end{remark}

We shall now show that Condition~$\textup{(Z}_0\textup{)}$ implies Condition~(W). We will need a few lemmas.

\begin{lemma}[see Part~I~{\cite[Corollary~4.1]{Po:Part1}}]\label{lem:Tauberian.Com-D2}
 For every $b\in \overline{\sA}$, we have $[D^{-2},b]\in (\sL_{\frac{p}{2},\infty})_0$.
\end{lemma}

\begin{lemma}[see Part~I~{\cite[Lemma~3.6]{Po:Part1}}]\label{lem:ST.ComparisonD2}
 For every $a\in \overline{\sA}_{++}$, we have 
\begin{equation*}
 \lim_{j\rightarrow \infty} j^{\frac{2}{p}}\lambda_{j}\left(a^{-1}D^{-2}a^{-1}\right) =  \lim_{j\rightarrow \infty} j^{\frac{2}{p}}\lambda_{j}\left(aD^{2}a\right)^{-1}, 
\end{equation*}
provided any of the above limits exists in $[0,\infty]$. That is, if one limit exists, then the other limit exists and the limits agree.  
 \end{lemma}

In what follows,  we set
\begin{equation}\label{eq:STC.sA++q}
 \sA_{++}^q:=\left\{a^{q}; \ a\in {\sA}_{++}\right\}, \qquad q\neq 0. 
\end{equation}
If $\sA$ is closed under holomorphic functional calculus, then $\sA_{++}^{q}=\sA_{++}$, but in general we do not have an equality for $q\neq 1$.

\begin{lemma}[see Part~I~{\cite[Lemma~3.7]{Po:Part1}}]\label{lem:STC.density}
If $a\in \overline{\sA}_+$, then, given any $q\neq 0$, we can always find a sequence $(a_\ell)_{\ell\geq 0}\subseteq \sA_{++}^q$ such that $(a_\ell)^r \rightarrow a^r$ in $\overline{\sA}$ for all $r>0$. In particular,  $\sA_{++}^q$ is a dense subset of $\overline{\sA}_+$. 
 \end{lemma}

We are now in a position to prove the main result of this subsection.

\begin{proposition}\label{prop:Tauberian.Weyl-W}
 Condition~$\textup{(Z}_0\textup{)}$ implies Condition~\textup{(W)}.
\end{proposition}
\begin{proof}
 Assume Condition~$\textup{(Z}_0\textup{)}$ holds, and set $\alpha=\min(1,p/2)$, so that $\alpha\in (0,1]$ and $\alpha<p$. By Lemma~\ref{lem:Tauberian.Z0-W} the Weyl laws~(\ref{eq:Tauberian.Weyl}) hold for all $a\in \sA_{++}$. 
 
 Let $a\in \sA_{++}$, and set $s=2\alpha^{-1}$, so that $s>1$ and $\alpha s=2$. Set $Y=a^{\frac12}|D|^{-\alpha}a^{\frac12}$. By Theorem~\ref{thm:Intro.CSZ} we have $Y^{s}-D^{-2}a^{s}\in \sL_{(3-\epsilon)^{-1}p,\infty}$ for all $\epsilon\in (0,1)$, and so $Y^{s}-D^{-2}a^{s}\in (\sL_{\frac{p}{2},\infty})_0$. Moreover, Lemma~\ref{lem:Tauberian.Com-D2} ensures that
 \begin{equation*}
 D^{-2}a^{s}-a^{\frac{s}{2}}D^{-2}a^{\frac{s}{2}}=\big[D^{-2},a^{\frac{s}{2}}\big]a^{\frac{s}{2}} \in (\sL_{\frac{p}{2},\infty})_0.
 \end{equation*}
 Thus, $Y^s-a^{\frac1\alpha}D^{-2}a^{\frac1\alpha}\in (\sL_{\frac{p}{2},\infty})_0$. Both operators are positive. Moreover, as $\lambda_j(Y^s)=\lambda_j(Y)^s$, the Weyl law~(\ref{eq:Tauberian.Weyl}) for $a$ gives
 \begin{equation*}
 \lim_{j\rightarrow \infty} j^{\frac{2}{p}}\lambda_j\big(Y^s\big)=\Big[\lim_{j\rightarrow \infty} j^{\frac{\alpha}{p}}\lambda_j(Y)\Big]^{\frac{2}{\alpha}}=\tau\big[a^{\frac{p}{\alpha}}\big]^{\frac{2}{p}}.
 \end{equation*}
 Proposition~\ref{prop:Bir-Sol} then ensures that
 \begin{equation}\label{eq:Tauberian.Weyl-cD-2c}
 \lim_{j\rightarrow \infty} j^{\frac{2}{p}}\lambda_j\big(cD^{-2}c\big)=\tau\big[c^{p}\big]^{\frac{2}{p}}
 \end{equation}
 for $c=a^{\frac1\alpha}$. That is, the Weyl laws~(\ref{eq:Tauberian.Weyl-cD-2c}) hold for all $c\in \sA_{++}^{1/\alpha}$.

 Let $c\in \overline{\sA}_+$. By Lemma~\ref{lem:STC.density} there is a sequence $(c_\ell)_{\ell\geq 0}\subseteq \sA_{++}^{1/\alpha}$ such that $c_\ell\rightarrow c$ and $(c_\ell)^p\rightarrow c^p$ in $\overline{\sA}$. Then $c_\ell D^{-2}c_\ell\rightarrow cD^{-2}c$ in $\sL_{\frac{p}{2},\infty}$. As $\tau$ is continuous, by using Proposition~\ref{prop:Bir-Sol} we get
 \begin{equation*}
 \lim_{j\rightarrow \infty} j^{\frac{2}{p}}\lambda_j\big(cD^{-2}c\big)=\lim_{\ell\rightarrow \infty}\tau\big[c_\ell^{p}\big]^{\frac{2}{p}}=\tau\big[c^{p}\big]^{\frac{2}{p}}.
 \end{equation*}
 Thus, the Weyl laws~(\ref{eq:Tauberian.Weyl-cD-2c}) hold for all $c\in \overline{\sA}_+$.

 Finally, let $a\in \sA_{++}$. Applying~(\ref{eq:Tauberian.Weyl-cD-2c}) to $c=a^{-1}$ and using Lemma~\ref{lem:ST.ComparisonD2} we get
 \begin{equation*}
 \lim_{j\rightarrow \infty} j^{\frac{2}{p}}\lambda_{j}\big(aD^{2}a\big)^{-1}=\lim_{j\rightarrow \infty} j^{\frac{2}{p}}\lambda_{j}\big(a^{-1}D^{-2}a^{-1}\big)=\tau\big[a^{-p}\big]^{\frac{2}{p}}.
 \end{equation*}
 As $\tau$ is a non-zero positive linear map and $a^{-p}$ is positive and invertible, we have $\tau[a^{-p}]>0$. This gives the Weyl law~(\ref{eq:Intro.conf-Weyl}) for $a$. That is, Condition~(W) is satisfied. The proof is complete.
\end{proof}

\subsection{Condition~(Z)} We may replace Condition ($\textup{Z}_0$) by a more standard Tauberian condition, which is easier to check in practice. Incidentally, Condition~(Z) is the type of Tauberian condition sought by Alain Connes. 

We have a continuous embedding $T\rightarrow \varphi_T$ of $\sL_1$ into $\overline{\sA}'$ given by
\begin{equation}\label{eq:Tauberian.phiT}
 \varphi_T(a)=\Tr\big[a T\big], \qquad a \in \overline{\sA}, \quad T\in \sL_1. 
\end{equation}
In particular, as $|D|^{-s}$, $\Re s>p$, is a holomorphic family of trace-class operators (\emph{cf}.\ Lemma~\ref{lem:Tauberian.powers}), we get a 
holomorphic family  $\varphi_s$, $\Re s>p$, in $\overline{\sA}'$ given by
\begin{equation*}
 \varphi_s(a):=\varphi_{|D|^{-s}}(a)=\Tr\big[a|D|^{-s}\big], \qquad a\in \overline{\sA}, \quad \Re s>p. 
\end{equation*}
If $\sB$ is any Fr\'echet subalgebra of $\overline{\sA}$ with strong dual $\sB'$, then, as we have a continuous embedding of $\overline{\sA}'$ into $\sB'$, we see that 
$(\varphi_s)_{|\sB}$, $\Re s>p$, is a holomorphic family in $\sB'$.

In what follows we denote by $\Omega_p$ the open halfplane $\{\Re z>p\}$ with closure $\overline{\Omega}_p:=\{\Re z\geq p\}$. 

\begin{condition}[Condition~\textup{(Z)}]
There is a holomorphically closed Fr\'echet subalgebra  $\sB\subseteq \overline{\sA}$ containing $\sA$ such that, for every $a\in \sB$,  
the function
 \begin{equation}\label{eq:Tauberian.phis}
\Omega_p\ni s \longrightarrow \varphi_s(a) -p\tau(a)(s-p)^{-1}
\end{equation}
 extends continuously to $\overline{\Omega}_p$. 
\end{condition}

\begin{remark}
 In practice we often may take $\sB=\overline{\sA}$, or $\sB=\sA$ if $\sA$ is itself a holomorphically closed Fr\'echet subalgebra. 
\end{remark}

\begin{remark}\label{rmk:Tauberian.(Z)}
 For $z\in \Omega_p$ and $a\in \sB$, set $\psi_z(a):=\varphi_z(a)-p\tau(a)(z-p)^{-1}$. Condition~(Z) means that, given any $a\in \sB$, for all $s\in \partial \Omega_p$, we have a limit
 \begin{equation*}
 \psi_s(a):=\lim_{\substack{z\rightarrow s\\ z\in \Omega_p}} \psi_z(a),
\end{equation*}
in such a way that the function $\overline{\Omega}_p\ni s\rightarrow \psi_s(a)$ is continuous. By Banach--Steinhaus theorem each $\psi_s$, $s\in \partial \Omega_p$, is a continuous linear functional on $\sB$. It follows that $(\psi_s)_{s\in \overline{\Omega}_p}$ is a continuous family in $\sB'_{\textup{weak}}$. 
\end{remark}

We have the following relationship between Condition~(Z) and Condition~$\textup{(Z}_0\textup{)}$. 

\begin{lemma}\label{lem:Tauberian.Z-Z0}
Condition~\textup{(Z)} implies Condition~$\textup{(Z}_0\textup{)}$.
\end{lemma}
\begin{proof}
As Condition~\textup{(Z)} holds, it follows from Remark~\ref{rmk:Tauberian.(Z)} that there is a holomorphically closed Fr\'echet subalgebra $\sB\subseteq \overline{\sA}$ and a continuous family $(\psi_s)_{s\in \overline{\Omega}_p}$ in $\sB'_{\textup{weak}}$ such that
\begin{equation}\label{eq:Tauberian.psis}
 \varphi_s(a)=p\tau(a)(s-p)^{-1} +\psi_s(a) \qquad \forall a\in \sB \quad \forall s\in \Omega_p. 
\end{equation}
 
 Let $a\in {\sB}_{++}$.  We have
 \begin{equation}\label{eq:Tauberian.as-ap}
 a^s-a^p=a^p\big(e^{(s-p)\log a}-1\big)=a^p \sum_{n\geq 1} \frac1{n!}(s-p)^n (\log a)^n=(s-p)b(s), 
\end{equation}
 where we have set
 \begin{equation*}
 b(s)=\sum_{n\geq 0} \frac{1}{(n+1)!} (s-p)^n a^p(\log a)^{n+1}, \qquad s\in \C. 
\end{equation*}
 In particular, $b(s)$, $s\in \C$, is an analytic family in $\overline{\sA}$. 
 
 Combining~(\ref{eq:Tauberian.psis}) with~(\ref{eq:Tauberian.as-ap}) shows that, for $\Re s>p$, we have
 \begin{align}
 \Tr\big[a^{s}|D|^{-s}\big]=\varphi_s(a^s) & = p\tau\big[a^{s}\big](s-p)^{-1} + \psi_s\big[a^{s}\big]\\
 &= p\tau\big[a^{p}\big](s-p)^{-1}+  p\tau\big[b(s)\big] +  \psi_s\big[a^{s}\big]. 
\end{align}
 
 The continuity of $\tau$ ensures that $\tau[b(s)]$ is an entire function on $\C$. Moreover, we know from Lemma~\ref{lem:Tauberian.hol-closed-powers} that $a^s$, $s\in \C$, is a holomorphic family in $\sB$. Let $s_0\in \overline{\Omega}_p$ and $K$ a compact neighborhood of $s_0$. As $(\psi_s)_{s\in\overline{\Omega}_p}$ is a continuous family in $\sB'_{\textup{weak}}$, its restriction $(\psi_s)_{s\in K}$ is weakly bounded, and hence is equicontinuous by the Banach--Steinhaus theorem. This means there is continuous seminorm $\fp$ on $\sB$ such that
\begin{equation*}
 |\psi_s(x)|\leq \fp(x) \qquad \forall x\in \sB \quad \forall s\in  K. 
\end{equation*}
Thus, for all $s\in K$, we have
\begin{align*}
 \left| \psi_s(a^s)-\psi_{s_0}(a^{s_0})\right| &\leq \left| \psi_s(a^s-a^{s_0})\right|+ \left| \psi_s(a^{s_0})-\psi_{s_0}(a^{s_0})\right| \\
 & \leq \fp(a^s-a^{s_0})+  \left| \psi_s(a^{s_0})-\psi_{s_0}(a^{s_0})\right|. 
\end{align*}
As $ \psi_s(a^{s_0})\rightarrow \psi_{s_0}(a^{s_0})$ and $a^{s}\rightarrow a^{s_0}$ in $\sB$ as $s\rightarrow s_0$, we then deduce that $\psi_s(a^s)\rightarrow \psi_{s_0}(a^{s_0})$ as $s\rightarrow s_0$. This shows that the function $s\rightarrow \psi_s(a^s)$ is continuous on $\overline{\Omega}_p$. 

It follows from all this that, for every $a\in \sB_{++}$, the function
 \begin{equation*}
\Omega_p\ni s \longrightarrow  \Tr\big[a^{s}|D|^{-s}\big] - p\tau\left[a^{p}\right](s-p)^{-1}
\end{equation*}
extends continuously to $\overline{\Omega}_p$.  That is, Condition~$\textup{(Z}_0\textup{)}$ is satisfied, proving the result.  
\end{proof}

\subsection{Condition~(H)} There is a close relationship between singularities of zeta functions and short-time behaviours of traces of heat semigroups (see, e.g., \cite{GS:JGA96}). In various examples the latter are often simpler to study. Therefore, it stands to reason to seek a heat semigroup version of Condition~$\textup{(Z)}$.

\begin{lemma}\label{lem:Tauberian.heat-semigroup}
 The heat semigroup $e^{-tD^2}$, $t>0$, yields a continuous family of trace-class operators such that
\begin{gather}\label{eq:Tauberian.heat-large-time}
 \big\| (1-\Pi_0)e^{-tD^2}\big\|_1 =\op{O}\big(e^{-t\lambda_0(D^2)}\big) \quad \text{as $t\rightarrow \infty$},\\
  \big\| e^{-tD^2}\big\|_1 =\op{O}\big(t^{-\frac{p}{2}}\big) \quad \text{as $t\rightarrow 0^+$}.
  \label{eq:Tauberian.heat-short-time}
\end{gather}
\end{lemma}
\begin{proof}
 This result is standard. We sketch a proof for the reader's convenience. Let $(\xi_j)_{j\geq 0}$ be an orthonormal family in $\sH$ such that $D^2\xi_j=\lambda_j(D^2)\xi_j$. For all $t>0$, we have
 \begin{equation}\label{eq:Tauberian.Schmidt-heat}
 e^{-tD^2} = \Pi_0 + \sum_{j\geq 0}e^{-t\lambda_j(D^2)} |{\xi_j}\rangle\langle \xi_j|, 
\end{equation}
where the series converges in $\sK$. As $D^{-1}\in \sL_{p,\infty}$, setting $C=\|D^{-1}\|_{p,\infty}^{-2}$ we get
\begin{equation*}
 \lambda_j\big(D^2\big)=\mu_j\big(D^{-1}\big)^{-2}\geq C(j+1)^{\frac{2}{p}} \quad \text{for all $j\geq 0$}. 
\end{equation*}
Thus, given any $b>0$, for all $j\geq 0$ and $t\geq b$, we have
\begin{equation}\label{eq:Tauberian.trace-norm-xij}
 \big\|e^{-t\lambda_j(D^2)} |{\xi_j}\rangle\langle \xi_j|\big\|_1 = e^{-t\lambda_j(D^2)} \leq e^{-bC(j+1)^{\frac{2}{p}}}. 
\end{equation}
This shows that the Schmidt series $\sum_{j\geq 0}e^{-t\lambda_j(D^2)} |{\xi_j}\rangle\langle \xi_j|$ converges normally in the Banach space $C_0([b,\infty); \sL_1)$ for all $b>0$. It then follows that 
$e^{-tD^2}$, $t>0$, is a continuous family in $\sL_1$. 

Note that $\|(1-\Pi_0)e^{-tD^2}\|=e^{-t\lambda_0(D^2)}$ for all $t>0$. Thus, for all $t>1$, by using the semigroup property $ e^{-tD^2}=e^{-(t-1)D^2}e^{-D^2}$, we get
\begin{align*}
  \big\| (1-\Pi_0)e^{-tD^2}\big\|_1 &= \big\| (1-\Pi_0)e^{-(t-1)D^2}e^{-D^2}\big\|_1 \\
  & \leq \big\| (1-\Pi_0)e^{-(t-1)D^2}\| \|e^{-D^2}\big\|_1\\
  & \leq  \|e^{-D^2}\big\|_1e^{-(t-1)\lambda_0(D^2)}. 
\end{align*}
This gives~(\ref{eq:Tauberian.heat-large-time}). 

Moreover, given any  $t>0$, by using~(\ref{eq:Tauberian.Schmidt-heat}) and~(\ref{eq:Tauberian.trace-norm-xij}) for $b=t$ we get
\begin{equation*}
 \big\| e^{-tD^2}\big\|_1 \leq \|\Pi_0\|_1 + \sum_{j\geq 0}  \big\|e^{-t\lambda_j(D^2)} |{\xi_j}\rangle\langle \xi_j|\big\|_1 \leq \dim \ker D + 
  \sum_{j\geq 0} e^{-tC(j+1)^{\frac{2}{p}}}.
\end{equation*}
Note that
\begin{equation*}
  \sum_{j\geq 0} e^{-tC(j+1)^{\frac{2}{p}}} \leq \int_0^\infty e^{-tCx^{\frac{2}{p}}}dx = (Ct)^{-\frac{p}{2}} \int_0^\infty e^{-x^{\frac{2}{p}}}dx.
\end{equation*}
This gives~(\ref{eq:Tauberian.heat-short-time}). The proof is complete. 
\end{proof}

As $e^{-tD^2}$, $t>0$, is a continuous family in $\sL_1$, it gives rise to a continuous family $\theta_t$, $t>0$, in $\overline{\sA}'$ given by
\begin{equation}\label{eq:Tauberian.thetat}
 \theta_t(a)= \Tr\big[ae^{-tD^2}\big], \qquad a \in \overline{\sA}, \quad t>0. 
\end{equation}
In particular, if $\sB$ is any Fr\'echet subalgebra of $\overline{\sA}$, then $(\theta_t)_{|\sB}$, $t>0$, is a continuous family in $\sB'$. 

The heat semigroup version of Condition~$\textup{(Z)}$ is provided by the following condition. 

\begin{condition}[Condition~$\textup{(H)}$]
There are $\delta>0$ and a holomorphically closed Fr\'echet subalgebra  $\sB\subseteq \overline{\sA}$ containing $\sA$ such that, for all $a\in \sB$, as $t\rightarrow 0^+$ we have 
\begin{equation*}
 \theta_t (a) = \Gamma\left(1+\frac{p}{2}\right) t^{-\frac{p}{2}}\left(\tau(a) + \op{O}\big(t^\delta\big)\right). 
\end{equation*}
\end{condition}

\begin{remark}\label{rmk:Tauberian.(H)}
 Equivalently, Condition~(H) means that the family of linear functionals, 
 \begin{equation*}
 \varepsilon_t(a):=t^{-\delta}\left[ \Gamma\left(1+\frac{p}{2}\right)^{-1}t^{\frac{p}{2}}\theta_t(a)- \tau(a)\right], \qquad a\in \sB, \quad t>0,
\end{equation*}
is weakly bounded in $\sB'$ on every interval $(0,c]$, $c>0$. Note also that $(\varepsilon_t)_{t>0}$ is a continuous family in $\sB'$. 
\end{remark}

\begin{lemma}\label{lem:Tauberian.H-Z}
 Condition~$\textup{(H)}$ implies Condition~\textup{(Z)}.
\end{lemma}
\begin{proof}
 Assume Condition~$\textup{(H)}$ holds. For $\Re s>0$, we have
 \begin{align*}
  \Gamma\left(\frac{s}{2}\right)|D|^{-s} =  \Gamma\left(\frac{s}{2}\right)\left(D^2\right)^{-\frac{s}{2}}
  &=  \int_0^\infty t^{\frac{s}{2}}\left(1-\Pi_0\right)e^{-tD^2}\frac{dt}{t}\\
  & =  \int_0^1t^{\frac{s}{2}}e^{-tD^2}\frac{dt}{t} - \frac{2}{s} \Pi_0 + R(s), 
  \end{align*}
 where we have set  
 \begin{equation*}
 R(s)= \int_1^\infty t^{\frac{s}{2}}\left(1-\Pi_0\right)e^{-tD^2}\frac{dt}{t}. 
\end{equation*}
We have also used the fact that
\begin{equation*}
  \int_0^1 t^{\frac{s}{2}}\Pi_0e^{-tD^2}\frac{dt}{t}= \int_0^1t^{\frac{s}{2}}\Pi_0\frac{dt}{t}=\frac{2}{s}\Pi_0. 
\end{equation*}

The above integrals converge in $\sL(\sH)$. However, for $\Re s>p$, Lemma~\ref{lem:Tauberian.heat-semigroup} ensures that $ t^{-1+s/2}\left(1-\Pi_0\right)e^{-tD^2}$, $t>0$, is an $L^1$-family in $\sL_1$, and so the integrals actually converge in $\sL_1$. Thus, for all $a\in \overline{\sA}$, we have
\begin{equation*}
  \Gamma\left(\frac{s}{2}\right)\Tr\big[a|D|^{-s}\big]  =  \int_0^1t^{\frac{s}{2}}\Tr\big[ ae^{-tD^2}\big] \frac{dt}{t} - \frac{2}{s} \Tr\big[ a\Pi_0\big] + \Tr\big[ a R(s)\big], \quad \Re s >p.   
\end{equation*}
Equivalently, with the notation of~(\ref{eq:Tauberian.phiT}),~(\ref{eq:Tauberian.phis}) and~(\ref{eq:Tauberian.thetat}) we have
\begin{equation}\label{eq:Tauberian.vphis-int01-Rs}
  \Gamma\left(\frac{s}{2}\right)\varphi_s(a)  =  \int_0^1t^{\frac{s}{2}}\theta_t(a) \frac{dt}{t} - \frac{2}{s} \Tr\big[ a\Pi_0\big] + \varphi_{R(s)}(a), \quad \Re s >p.   
\end{equation}

Note that $(s,t)\rightarrow t^{-1+s/2}(1-\Pi_0)e^{-tD^2}$ is a continuous map from $\C\times [1,\infty)$ to $\sL_1$ which is holomorphic with respect to $s$. Moreover, thanks to~(\ref{eq:Tauberian.heat-large-time}), for any compact $K\subseteq \C$, there are $\beta_K\in \R$ and $C_K>0$ such that, for all $s\in K$ and $t\geq 1$, we have
\begin{equation*}
\big\|t^{-1+s/2}(1-\Pi_0)e^{-tD^2}\big\|_1 \leq  t^{\beta_K} \big\|(1-\Pi_0)e^{-tD^2}\big\|_1\leq C_K t^{\beta_K}e^{-t\lambda_0(D^2)}. 
\end{equation*}
As $t^{\beta_K}e^{-t\lambda_0(D^2)} \in L^1[1,\infty)$, it follows that $R(s)$, $s\in \C$,  is a holomorphic family in $\sL_1$, and so 
$\varphi_{R(s)}$, $s\in \C$, is a holomorphic family in $\overline{\sA}'$. 

As Condition~(H) holds, we know from Remark~\ref{rmk:Tauberian.(H)}, there are $\delta>0$, a holomorphically closed Fr\'echet subalgebra $\sB\subseteq \overline{\sA}$, and a weakly bounded continuous family $(\varepsilon_{t})_{0<t\leq 1}$ in $\sB'$ such that, for all $a\in \sB$, we have
\begin{equation*}
 \theta_t(a) = \Gamma\left(1+\frac{p}{2}\right)t^{-\frac{p}{2}}\left[\tau(a)  + t^{\delta}\varepsilon_t(a)\right], \qquad 0<t\leq 1.  
\end{equation*}
Thus, if $a\in \sB$, then, for $\Re s >p$, we have
\begin{align}
\int_0^1t^{\frac{s}{2}}\theta_t(a) \frac{dt}{t} 
 &=  \Gamma\left(1+\frac{p}{2}\right)\int_0^1 t^{\frac{s-p}{2}}\tau(a)  \frac{dt}{t} +  \Gamma\left(1+\frac{p}{2}\right)\int_0^1t^{\frac{s-p}{2}}t^{\delta-1}\varepsilon_t(a) dt \nonumber\\
 &= p\Gamma\left(\frac{p}{2}\right)\tau(a)(s-p)^{-1}+ \hat{\varepsilon}_s(a), 
 \label{eq:Tauberian.int01thetat}
 \end{align}
where we have set $\hat{\varepsilon}_s(a):=\Gamma(1+\frac{p}{2})\int_0^1 t^{(s-p)/2}t^{\delta-1}\varepsilon_t (a)dt$. As $\varepsilon_t(a)$ is bounded on $(0,1]$, the above integral converges for $\Re s>p-2\delta$, and so we see that 
$\hat{\varepsilon}_s(a)$ is a holomorphic function on the half-plane $\Re s>p-\delta$.

Combining~(\ref{eq:Tauberian.vphis-int01-Rs}) and~(\ref{eq:Tauberian.int01thetat}) shows that, for all $a\in \sB$, we have
\begin{equation*}
 \Gamma\left(\frac{s}{2}\right)\varphi_s(a) = p\Gamma\left(\frac{p}{2}\right)\tau(a) (s-p)^{-1}+ \psi_s(a), \qquad \Re s>p. 
\end{equation*}
where we have set
\begin{equation*}
 \psi_s(a):= \hat{\varepsilon}_s(a) - \frac{2}{s} \Tr\big[ a\Pi_0\big] + \varphi_{R(s)}(a), \qquad  \Re s>p-\delta. 
\end{equation*}
Note that $s\rightarrow \psi_s(a)$ is a holomorphic function on the halfplane $\Re s> p-\delta$.

In addition, as $\Gamma(s/2)^{-1}$ is a holomorphic function on the half-plane $\Re s>0$, we may write
\begin{equation*}
 \Gamma\left(\frac{s}{2}\right)^{-1} =\Gamma\left(\frac{p}{2}\right)^{-1}+(s-p)u(s), 
\end{equation*}
where $u(s)$ is a holomorphic function on the half-plane $\Re s>0$. Thus, for $\Re s>p$, we have
\begin{align*}
\varphi_s(a) & = p\Gamma\left(\frac{s}{2}\right)^{-1}\Gamma\left(\frac{p}{2}\right)\tau(a) (s-p)^{-1}+ 
 \Gamma\left(\frac{s}{2}\right)^{-1} \psi_s(a)\\
 &= p\tau(a) (s-p)^{-1} + pu(s)\Gamma\left(\frac{p}{2}\right)\tau(a) +  \Gamma\left(\frac{s}{2}\right)^{-1} \psi_s(a). 
\end{align*}
Equivalently,
\begin{equation*}
 \varphi_s(a)-p\tau(a) (s-p)^{-1} =  pu(s)\Gamma\left(\frac{p}{2}\right)\tau(a) +  \Gamma\left(\frac{s}{2}\right)^{-1} \psi_s(a). 
\end{equation*}
As the r.h.s.\ is a holomorphic function on the halfplane $\Re s>p-\delta$, this shows that the l.h.s.\ has a holomorphic extension to that halfplane for all $a\in \sB$. This ensures that Condition~$\textup{(Z)}$ holds, proving the result. 
\end{proof}

\begin{remark}
 The proof above shows that Condition~(H) implies a slightly stronger form of Condition~(Z), since we obtain an analytic continuation of the functions~(\ref{eq:Tauberian.phis}) to the halfplane $\{\Re z>p-\delta\}$.  
\end{remark}

\subsection{Proof of the Tauberian Theorem}
 By Lemma~\ref{lem:Tauberian.H-Z} Condition~(H) implies Condition~(Z), and by Lemma~\ref{lem:Tauberian.Z-Z0} Condition~(Z) implies Condition $\textup{(Z}_0\textup{)}$. By Proposition~\ref{prop:Tauberian.Weyl-W} the latter implies Condition~(W). Therefore, Condition~(Z) and Condition~(H) both imply Condition~(W). This completes the proof of Theorem~\ref{thm:Intro.Tauberian}.\hfill$\Box$

\appendix 
\section{Closedness Under Holomorphic Functional Calculus}\label{app:hol-extension}
Let $(\sA,\sH,D)$ be a spectral triple, and denote by $\overline{\sA}$ the closure of $\sA$ in $\sL(\sH)$. The aim of this appendix is to prove the following result. 

\begin{proposition}\label{prop:App.hol-extension}
 There is a $*$-subalgebra $\sA_1$ of $\overline{\sA}$ containing $\sA$ such that $(\sA_1,\sH,D)$ is a spectral triple and $\sA_1$ is closed under holomorphic functional calculus.
\end{proposition}
\begin{proof}
 Let $\sA_1$ be the set of operators $T\in \overline{\sA}$ such that $T(\dom(D))\subseteq \dom(D)$ and $[D,T]\in \sL(\sH)$. As $(\sA,\sH,D)$ is a spectral triple, $\sA_1$ contains $\sA$. 

If $S,T\in \sA_1$, then $ST(\dom(D))\subseteq S(\dom(D))\subseteq \dom(D)$, and on $\dom(D)$ we have $[D,ST]=S[D,T]+[D,S]T$. Thus, $ST\in \sA_1$, and so $\sA_1$ is a subalgebra of $\overline{\sA}$. Let $T\in \sA_1$. If $\xi,\eta\in \dom(D)$, then
\begin{equation*}
 \scal{T^*\xi}{D\eta}=\scal{\xi}{DT\eta}-\scal{\xi}{[D,T]\eta}=\scal{T^*D\xi-[D,T]^*\xi}{\eta}. 
\end{equation*}
This shows that $T^*\xi \in \dom(D^*)=\dom(D)$ and $DT^*\xi=T^*D\xi-[D,T]^*\xi$. Thus, $T^*(\dom(D))\subseteq \dom(D)$ and $[D,T^*]=-[D,T]^*$ on $\dom(D)$. As $T^*\in \overline{\sA}$, we see that $T^*\in \sA_1$, and so $\sA_1$ is a $*$-subalgebra of $\overline{\sA}$. It then follows that $(\sA_1,\sH,D)$ is a spectral triple. 

We equip $\sA_1$ with the norm, 
\begin{equation*}
 \|T\|_1:= \|T\| + \left\|[D,T]\right\|, \qquad T\in \sA_1. 
\end{equation*}
The Leibniz rule $[D,ST]=S[D,T]+[D,S]T$ gives $\|ST\|_1\leq \|S\|_1\|T\|_1$ for all $S,T\in \sA_1$. Let $(T_\ell)_{\ell\geq 0}$ be a Cauchy sequence in $\sA_1$. Then $(T_\ell)_{\ell\geq 0}$ and $([D,T_\ell])_{\ell\geq 0}$ are Cauchy sequences in $\sL(\sH)$, and so there are $T\in \overline{\sA}$ and $S\in \sL(\sH)$ such that $T_\ell \rightarrow T$ and $[D,T_\ell]\rightarrow S$ in $\sL(\sH)$. Let $\xi \in \dom(D)$. Then $T_\ell \xi \rightarrow T\xi$, and 
\begin{equation*}
 D\left(T_\ell \xi\right)=T_\ell(D\xi)+ \big[D,T_\ell\big]\xi \longrightarrow T(D\xi)+S\xi. 
\end{equation*}
As $D$ is selfadjoint, its graph is closed, and so $T\xi \in \dom(D)$ and $D(T\xi)=T(D\xi)+S\xi$, i.e., $[D,T]\xi=S\xi$. Thus, $T\in \sA_1$ and $[D,T]=S$, and so $T_\ell\rightarrow T$ in $\sA_1$. Therefore, $\sA_1$ is a Banach $*$-algebra. As $\|T\|\leq \|T\|_1$, it is a Fr\'echet subalgebra of $\overline{\sA}$ in the sense of Section~\ref{subsec:Prelim.ST}. 

It remains to show that $\sA_1$ is closed under holomorphic functional calculus. As mentioned in Section~\ref{subsec:Prelim.ST}, it is enough to show that $\sA_1$ is inverse closed. This is a consequence of well-known results of Blackadar--Cuntz~\cite[Proposition~3.12, Example~3.16(a)]{BC:JOT91} and 
Connes~\cite[Lemma~III.6.2]{Co:NCG} on the holomorphic closedness of domains of closed derivations. We actually have a slightly more general result. 

\begin{lemma}\label{lem:App.inverse-closed}
 Let $T\in \sL(\sH)$ be such that $T(\dom(D))\subseteq \dom(D)$ and $[D,T]\in \sL(\sH)$. If $T$ is invertible in $\sL(\sH)$, then $T^{-1}(\dom(D))\subseteq \dom(D)$, and on $\dom(D)$ we have
 \begin{equation}\label{eq:App.com-inverse}
  [D,T^{-1}]=-T^{-1}[D,T]T^{-1}.
 \end{equation}
 In particular, $[D,T^{-1}]\in \sL(\sH)$.
\end{lemma}

\begin{proof}[Proof\protect\footnotemark]\footnotetext{The proof of this lemma was obtained with the help of Claude AI (Opus~5.5).}
 Set $S:=T^{-1}$. Given $\mu\in \R\setminus\{0\}$, set $R_\mu:=(D-i\mu)^{-1}$. As $D$ is selfadjoint, $R_\mu$ is a bounded operator with $\|R_\mu\|\leq |\mu|^{-1}$ and $R_\mu(\sH)=\dom(D)$. Moreover, $(D-i\mu)R_\mu=1$ on $\sH$ and $R_\mu(D-i\mu)=1$ on $\dom(D)$. As $T(\dom(D))\subseteq \dom(D)$, the operator $TR_\mu$ maps $\sH$ to $\dom(D)$, and so
 \begin{equation*}
  R_\mu T-TR_\mu = R_\mu T(D-i\mu)R_\mu- R_\mu(D-i\mu)TR_\mu=-R_\mu[D,T]R_\mu .
 \end{equation*}
 Multiplying this equality on both sides by $S$ gives $SR_\mu-R_\mu S=-SR_\mu[D,T]R_\mu S$, that is,
 \begin{equation*}
  SR_\mu\big(1+[D,T]R_\mu S\big)=R_\mu S.
 \end{equation*}
 Choose $\mu$ so that $|\mu|>\|[D,T]\|\|S\|$. Then $\|[D,T]R_\mu S\|<1$, and so $1+[D,T]R_\mu S$ is invertible in $\sL(\sH)$. Therefore,
 \begin{equation*}
  SR_\mu =R_\mu S\big(1+[D,T]R_\mu S\big)^{-1}.
 \end{equation*}
 As $R_\mu(\sH)=\dom(D)$, we deduce that $S(\dom(D))=SR_\mu(\sH)\subseteq R_\mu(\sH)=\dom(D)$.

 Let $\xi\in \dom(D)$. Then $\eta:=S\xi$ is in $\dom(D)$, and so $[D,T]\eta=DT\eta-TD\eta=D\xi-TDS\xi$. Applying $S$ gives $S[D,T]S\xi=SD\xi-DS\xi$, i.e., $[D,S]\xi=-S[D,T]S\xi$. This proves~(\ref{eq:App.com-inverse}), and the boundedness of $[D,S]$ follows. The proof of the lemma is complete.
\end{proof}

Let $T\in \sA_1$ be invertible in $\overline{\sA}$. Then $T^{-1}\in \overline{\sA}$, and Lemma~\ref{lem:App.inverse-closed} ensures that $T^{-1}(\dom(D))\subseteq \dom(D)$ and $[D,T^{-1}]\in \sL(\sH)$. Thus, $T^{-1}\in \sA_1$. This shows that $\sA_1$ is inverse closed, and hence it is closed under holomorphic functional calculus. The proof of Proposition~\ref{prop:App.hol-extension} is complete. 
\end{proof}

\begin{remark}\label{rmk:App.hol-extension}
 The proof of Proposition~\ref{prop:App.hol-extension} shows that the algebra $\sA_1$ is a Banach $*$-algebra with respect to the norm $\|T\|_1=\|T\|+\|[D,T]\|$. Therefore, it is a holomorphically closed Fr\'echet subalgebra of $\overline{\sA}$ containing $\sA$, and $T\rightarrow [D,T]$ is a continuous linear map from $\sA_1$ to $\sL(\sH)$. 
\end{remark}

\begin{remark}\label{rmk:App.hol-extension-B}
 If $\sB\subseteq \overline{\sA}$ is a holomorphically closed Fr\'echet subalgebra containing $\sA$, then $\sB\cap \sA_1$ is a holomorphically closed Fr\'echet subalgebra of $\overline{\sA}$ containing $\sA$ with respect to the topology defined by the seminorms of $\sB$ and the norm $\|\cdot\|_1$. Indeed, a Cauchy sequence in $\sB\cap \sA_1$ converges in $\sB$ and in $\sA_1$, and the two limits agree, since both topologies are stronger than that of $\overline{\sA}$. Moreover, if $T\in \sB\cap \sA_1$ is invertible in $\overline{\sA}$, then $T^{-1}$ is contained in $\sB$ and in $\sA_1$. 
\end{remark}

\begin{remark}
 The algebra $\sA_1$ is the largest $*$-subalgebra of $\overline{\sA}$ satisfying Proposition~\ref{prop:App.hol-extension}. In the case of a Dirac spectral triple $(C^\infty(M), L^2(M,\shS),\shD)$, the algebra $\sA_1$ is the algebra of Lipschitz functions on $M$. 
 For other applications (e.g., index theory) it may be convenient to have a $*$-subalgebra $\sA_1$ which is as small as possible. One way to get such an algebra is to consider the intersection of the domains of all ``higher'' derivations containing $\sA$, as in~\cite{Mo:JFA24}.   
\end{remark}

\begin{remark}
 A version of Proposition~\ref{prop:App.hol-extension} for smooth (possibly nonunital) spectral triples was obtained by Rennie~\cite[Lemma~16]{Re:KT03}. Our approach does not require any smoothness assumption. 
\end{remark}

\section{Proof of Proposition~\ref{lem:Com.estimate-Delta-weak}}\label{app:BS}
In this appendix, for the reader's convenience, we include a self-contained proof of Proposition~\ref{lem:Com.estimate-Delta-weak}. 
We use the notation of Section~\ref{chap:DOIs}. 

\begin{proof}[Proof of Proposition~\ref{lem:Com.estimate-Delta-weak}] 
 The proof follows the outline of the proof of~\cite[Theorem~4.1]{BS:IEOT03} in our special setting. We first show that there is $C_{\delta \delta'}>0$ such that
 \begin{equation}\label{eq:Com.estimate-Delta2}
 \int_0^\infty\mu^{1+\alpha} \left|\bigscal{\Delta_{\delta,\delta'}(T;\mu)\xi}{\eta}\right| d\mu \leq C_{\delta \delta'}\|T\|\|\xi\|\|\eta\| \qquad \forall T\in \sL(\sH), \ \forall \xi,\eta\in \sH.  
\end{equation}
Let $T\in \sL(\sH)$ and $\xi,\eta\in \sH$. We have
\begin{equation*}
 \left| \bigscal{\Delta_{\delta,\delta'}(T;\mu)\xi}{\eta}\right| = \left|\bigscal{TR_{\delta'}(\mu)\xi}{R_\delta(\mu)\eta}\right|\leq 
 \|T\| \|R_{\delta'}(\mu)\xi\|\|R_{\delta}(\mu)\eta\|. 
\end{equation*}
 Moreover, as  $\alpha+\delta+\delta'=2$, we have $1+\alpha=3-(\delta+\delta')$. Thus, 
 \begin{multline}\label{eq:Com.estimate-int-Delta-xieta}
 \int_0^\infty \mu^{1+\alpha} \left|\bigscal{\Delta_{\delta,\delta'}(T;\mu)\xi}{\eta}\right| d\mu 
 \leq \|T\| \int_0^\infty \mu^{3-(\delta+\delta')} \|R_{\delta'}(\mu)\xi\|\|R_{\delta}(\mu)\eta\| d\mu  \\
 \leq \|T\| \left(\int_0^\infty \mu^{3-2\delta'} \|R_{\delta'}(\mu)\xi\|^2 d\mu\right)^{\frac12} 
 \left(\int_0^\infty \mu^{3-2\delta} \|R_{\delta}(\mu)\eta\|^2 d\mu\right)^{\frac12}.  
\end{multline}
 
Bearing this in mind, set $c=\min\{|\lambda|; \ \lambda \in \Sp(D)\setminus 0\}$. Note that $c>0$. Denote by $\chi(\lambda)$ the characteristic function of 
$(-\infty,-c]\cup[c,\infty)$, and, for $0\leq \delta \leq 2$, set
\begin{equation*}
 r_\delta(\lambda;\mu)=\chi(\lambda)|\lambda|^\delta \left(\lambda^2+\mu^2\right)^{-1}, \qquad \lambda\in \R, \ \mu\in \R. 
\end{equation*}
 We then have $R_\delta(\mu)=r_\delta(D;\mu)$, i.e., 
 \begin{equation*}
 R_\delta(\mu) = \int_\R r_\delta(\lambda;\mu)dE(\lambda), \qquad \mu\in \R,
\end{equation*}
 where $E(\lambda)$ is the spectral measure of $D$. Therefore, given any $\xi\in \sH$, we have 
\begin{equation*}
 \left\|R_\delta(\mu)\xi\right\|^2=\scal{R_\delta(\mu)^2\xi}{\xi}=\int r_\delta(\lambda;\mu)^2 d\scal{E(\lambda)\xi}{\xi}= \int  
 \chi(\lambda)|\lambda|^{2\delta} \left(\lambda^2+\mu^2\right)^{-2}d\scal{E(\lambda)\xi}{\xi}.  
\end{equation*}
Thus, 
\begin{align*}
 \int_0^\infty \mu^{3-2\delta} \|R_\delta(\mu)\xi\|^2 d\mu&  = 
\int_0^\infty \int_\R  \mu^{3-2\delta}\chi(\lambda)|\lambda|^{2\delta} \left(\lambda^2+\mu^2\right)^{-2}d\scal{E(\lambda)\xi}{\xi}d\mu\\ 
 & =\int_\R \chi(\lambda) |\lambda|^{2\delta} \left(\int_0^\infty  \mu^{3-2\delta}\left(\lambda^2+\mu^2\right)^{-2}d\mu\right)d\scal{E(\lambda)\xi}{\xi}. 
\end{align*}
Note that, for $\lambda\neq 0$, we have  
\begin{equation*}
\int_0^\infty  \mu^{3-2\delta}\left(\lambda^2+\mu^2\right)^{-2}d\mu = |\lambda|^{-2\delta} \int_0^\infty \mu^{3-2\delta}\left(1+\mu^2\right)^{-2}d\mu. 
\end{equation*}
Therefore, we get 
\begin{equation}\label{eq:Com.int-Rdelta2}
 \int_0^\infty \mu^{3-2\delta} \|R_\delta(\mu)\xi\|^2 d\mu
 = \int_0^\infty  \mu^{3-2\delta}\left(1+\mu^2\right)^{-2}d\mu \int_\R \chi(\lambda)d\scal{E(\lambda)\xi}{\xi}. 
\end{equation}

We observe that
\begin{equation}\label{eq:Com.int-chi}
 \int_\R \chi(\lambda)d\scal{E(\lambda)\xi}{\xi}=\scal{\chi(D)\xi}{\xi}=\scal{(1-\Pi_0)\xi}{\xi}\leq \|\xi\|^2.
\end{equation}
Moreover, the integral $I(\delta):= \int_0^\infty  \mu^{3-2\delta}\left(1+\mu^2\right)^{-2}d\mu$ is finite if and only if the following two conditions are satisfied:
\begin{itemize}
 \item $3-2\delta>-1$, i.e., $\delta<2$. 
 
 \item $(3-2\delta)-4=-1-2\delta<-1$, i.e., $\delta>0$. 
\end{itemize}
 Therefore, if $\delta\in (0,2)$, then $I(\delta)<\infty$, and so by using~(\ref{eq:Com.int-Rdelta2}) and~(\ref{eq:Com.int-chi}) we get
 \begin{equation}\label{eq:Com.estimate-int-muRdxi}
 \int_0^\infty \mu^{3-2\delta} \|R_\delta(\mu)\xi\|^2 d\mu \leq I(\delta) \|\xi\|^2<\infty. 
\end{equation}
By~(\ref{eq:Com.estimate-int-Delta-xieta}) and~(\ref{eq:Com.estimate-int-muRdxi}) we then get
\begin{equation*}
  \int_0^\infty \mu^{1+\alpha} \left| \bigscal{\Delta_{\delta,\delta'}(T;\mu)\xi}{\eta}\right| d\mu \leq \sqrt{I(\delta)I(\delta')}\|T\|\|\xi\|\|\eta\|<\infty.  
\end{equation*}
This yields the estimate~(\ref{eq:Com.estimate-Delta2}). 

It follows from~(\ref{eq:Com.estimate-Delta2}) that the right-hand side of~(\ref{eq:Com.def-Psidd'-weak}) is well defined and is a sesquilinear form in $(\xi,\eta)$ bounded by $C_{\delta\delta'}\|T\|\|\xi\|\|\eta\|$. By the Riesz representation theorem there is a unique $\Psi^\alpha_{\delta,\delta'}(T)\in\sL(\sH)$ satisfying~(\ref{eq:Com.def-Psidd'-weak}), and $\|\Psi^\alpha_{\delta,\delta'}(T)\|\leq C_{\delta\delta'}\|T\|$. Its linearity in $T$ is immediate. The proof is complete. 
\end{proof}


\begin{thebibliography}{99}


\bibitem{BJ:CRAS83} Baaj, S.; Julg, P.: \emph{Th\'eorie bivariante de Kasparov et op\'erateurs non born\'es dans les $C^*$-modules hilbertiens}. 
C.\ R.\ Acad.\ Sci.\ Paris S\'er.\ I Math.\ \textbf{296} (1983), 875--878.
 


\bibitem{BS:PMF66} Birman, M.Sh.; Solomyak, M.Z.: \emph{Double Stieltjes operator integrals} (Russian). Problemy Mat.\ Fiz.\ No.\ 1, Izdat.\ Leningrad Univ., 1966, pp.\ 33--67.  
  
\bibitem{BS:FAA70} Birman, M.Sh.; Solomyak, M.Z.: \emph{The leading term of the spectral asymptotics for ``non-smooth'' elliptic problems}.  Funct.\ Anal.\ Appl.\ \textbf{4} (1970), 265--275.


\bibitem{BS:PMF73} Birman, M.Sh.; Solomyak, M.Z.: \emph{Double Stieltjes operator integrals. III} (Russian). Problemy Mat.\ Fiz.\ No.\ 6, Izdat.\ Leningrad Univ., 1973, pp.\ 27--53.


\bibitem{BS:Book} Birman, M.Sh.; Solomyak, M.Z.: \emph{Spectral theory of self-adjoint operators in Hilbert space}, D.\ Reidel, Boston, MA, 1987.


\bibitem{BS:JSM92} Birman, M.Sh.; Solomyak, M.Z.: \emph{Estimates for the difference of fractional powers of selfadjoint operators under unbounded perturbations}. 
J.\ Sov.\ Math.\ \textbf{61}, no.\ 2 (1992), 2018--2035.


\bibitem{BS:IEOT03}  Birman, M.Sh.; Solomyak, M.Z.: \emph{Double operator integrals in a Hilbert space}. Int.\ Eq.\ Oper.\ Th.\ \textbf{47} (2) (2003), 131--168.

\bibitem{BC:JOT91} Blackadar, B.; Cuntz, J.: \emph{Differential Banach algebra norms and smooth subalgebras of $C^*$-algebras}. J.\ Operator Theory \textbf{26} (1991), 255--282.

\bibitem{CP:CJM98} Carey, A.; Phillips, J.:  \emph{Unbounded Fredholm modules and spectral flow}. Can.\ J.\ Math.\ \textbf{50}(4) (1998), 673--718.


\bibitem{Co:NCG} Connes, A.: \emph{Noncommutative geometry}. Academic Press,\ Inc.,\ San Diego, CA, 1994.


\bibitem{CM:CRAS86} Connes, A.; Moscovici, H.: \emph{Transgression du caract\`ere de Chern et cohomologie cyclique}. 
C.\ R.\ Acad.\ Sci.\ Paris S\'er.\ I Math.\ \textbf{303} (1986), 913--918.


\bibitem{CSZ:MS17} Connes, A.;  Sukochev, F.;  Zanin, D.: \emph{Trace theorem for quasi-Fuchsian groups}. Mat.\ Sb.,
\textbf{208} (2017), 59--90.


\bibitem{Fa:PNAS51} Fan, K.: \emph{Maximum properties and inequalities for the eigenvalues of completely continuous operators}. Proc.\ Nat.\ Acad.\ Sci.\ U.S.A.\ \textbf{37} (1951), 760--766. 


\bibitem{GK:AMS69} Gohberg, I.C.; Krein, M.G.:
\emph{Introduction to the theory of linear nonselfadjoint operators}. Translations of Mathematical Monographs, Vol.\ 18. American Mathematical Society, Providence, R.I., 1969.


\bibitem{GS:JGA96} Grubb, G.; Seeley, R.T.: \emph{Zeta and eta functions for Atiyah--Patodi--Singer operators}. J.\ Geom.\ Anal.\  \textbf{6}(1) (1996),  31--77.

 


\bibitem{HP:AMS74} Hille, E.; Phillips, R.: \emph{Functional analysis and semigroups}, 3rd edition. Amer.\ Math.\ Soc., Providence, R.I., 1974.


 
 
  


\bibitem{MSZ:LMP22} McDonald, E.; Sukochev, F.; Zanin, D.: \emph{Semiclassical Weyl law and exact spectral asymptotics in noncommutative geometry}. 
Lett.\ Math.\ Phys.\ \textbf{112} (2022), Paper No.~77, 45 pp. 


\bibitem{Mo:JFA24}  Mohsen, O.: \emph{Higher order derivations on $C^*$-algebras and applications to smooth functional calculus and Schwartz functions on the tangent groupoid}. J.\ Funct.\ Anal.\ \textbf{287} (2024), no.~11, Paper No.~110615.


\bibitem{Pe:FAA85}  Peller, V.V.: \emph{Hankel operators in the perturbation theory of unitary and self-adjoint operators}. Funct.\ Anal.\ Appl.\ \textbf{19} (1985), 111--123.

\bibitem{Pf:PAMS} Pfister, H.: \emph{Continuity of the inverse}. Proc.\ Amer.\ Math.\ Soc.\ \textbf{95} (1985), 312--314. 


\bibitem{Po:Part1} Ponge, R.: \emph{Noncommutative geometry, spectral asymptotics, and semiclassical analysis, I. Main results}. 
Preprint \texttt{arXiv:2604.15008}. 


\bibitem{Re:KT03} Rennie, A.: \emph{Smoothness and locality for nonunital spectral triples}. K-Theory \textbf{28} (2003), 127--165.

\bibitem{SWW:CRAS98} Schrohe, E.; Walze, M.; Warzecha, J.-M.:
\emph{Construction de triplets spectraux \`a partir de modules de Fredholm}. 
C.\ R.\ Acad.\ Sci.\ Paris S\'er.\ I Math.\ \textbf{326} (1998), 1195--1199.


\bibitem{Sh:Springer01} Shubin, M.A.: \emph{Pseudodifferential operators and spectral theory}, 2nd edition. Springer, Berlin, 2001.


\bibitem{Si:AMS05} Simon, B.: \emph{Trace ideals and their applications}. Second edition. 
Mathematical Surveys and Monographs, 120. American Mathematical Society, Providence, RI, 2005.


 \bibitem{ST:LNM19} Skripka, A.; Tomskova, A.: \emph{Multilinear Operator Integrals: Theory and Applications}. 
 Lecture Notes in Mathematics,  2250, Springer, 2019, xi+192 pp.


\bibitem{SZ:PAMS21} Sukochev, F.; Zanin, D.: \emph{Optimal constants in non-commutative H\"older inequality for quasi-norms}. Proc.\ Amer.\ Math.\ Soc.\ \textbf{149} (2021), no.\ 9, 3813--3817.

 

\bibitem{SZ:Ast23} Sukochev, F.; Zanin, D.: \emph{The Connes character formula for locally compact spectral triples}. Ast\'erisque, no.\ 445, viii+150 pp.\ (2023).


\bibitem{Ze:BAPS} \.{Z}elasko, W.: \emph{A theorem on $B_0$ division algebras}. Bull.\ Acad.\ Polon.\ Sci.\ \textbf{8} (1960), 373--375.

\end{thebibliography}
\end{document}